\documentclass[leqno,12pt]{article}

\usepackage[dvips]{graphics}
\usepackage{color}
\usepackage{ifthen}
\usepackage{soul}
\usepackage{latexsym}
\usepackage{amsthm,amssymb}
\usepackage{amsbsy,amsfonts,amsmath}

\usepackage{setspace}

\newcommand{\bysame}{\leavevmode\hbox
to 3em{\hrulefill}\,}

\usepackage{a4}

\numberwithin{equation}{section}

\newtheorem{theorem}{Theorem}[section]
\newtheorem{proposition}[theorem]{Proposition}

\newtheorem{lemma}[theorem]{Lemma}

\theoremstyle{definition}

\theoremstyle{remark}
\newtheorem{remark}[theorem]{Remark}

\newcommand{\al}{\alpha}
\newcommand{\bN}{{\mathbb{N}}}
\newcommand{\bZ}{{\mathbb{Z}}}
\newcommand{\bZgeqo}{\bZ_{\geq 0}}
\newcommand{\bZleqo}{\bZ_{\leq 0}}
\newcommand{\bQ}{{\mathbb{Q}}}
\newcommand{\bR}{{\mathbb{R}}}
\newcommand{\bC}{{\mathbb{C}}}
\newcommand{\bCt}{\bC^\times}
\newcommand{\bCtinf}{\bCt_\infty}

\newcommand{\bK}{{\mathbb{K}}}
\newcommand{\bKt}{{\bK^\times}}
\newcommand{\bKtinf}{{\bK^\times_\infty}}
\newcommand{\Typrmid}{{\mathrm{id}}}
\newcommand{\TyprmSpan}{{\mathrm{Span}}}
\newcommand{\TyprmChar}{{\mathrm{Char}}}

\newcommand{\TypIntQ}{Q}
\newcommand{\TypIntacuteQ}{{\acute{\TypIntQ}}}

\newcommand{\TypIntKhR}{R}

\newcommand{\TypIntR}{{\mathcal{R}}}
\newcommand{\TypIntRp}{\TypIntR^+}
\newcommand{\TypIntRpreal}{\TypIntRp_{{\mathrm{real}}}}
\newcommand{\TypIntRpnull}{\TypIntRp_{{\mathrm{null}}}}

\newcommand{\TypIntmcMOm}{{\mathcal{M}}(\Omega)}
\newcommand{\TypIntmcLOm}{{\mathcal{L}}(\Omega)}

\newcommand{\Typkpch}{{\hat o}}

\newcommand{\Typpial}{\al}

\newcommand{\TypTHREEdummies}[3]{{#1(#2,#3)}}
\newcommand{\TypSecondTHREEdummies}[3]{#1^{#2,#3}}
\newcommand{\TypFOURdummies}[4]{{#1^{#2,#3}_{#4}}}
\newcommand{\TypSecondFOURdummies}[4]{{#1^{#2,#3}(#4)}}

\newcommand{\Typchi}{\chi}
\newcommand{\Typl}{\theta}
\newcommand{\Typlp}{\Typl^\prime}
\newcommand{\TypfkI}{I}
\newcommand{\TypfkJ}{J}
\newcommand{\Typpi}{\pi}

\newcommand{\TypfkAoriginal}{{\mathfrak{A}}}
\newcommand{\dummyTypfkAoriginal}[1]{\TypfkAoriginal_{#1}}
\newcommand{\TypfkAoriginalpi}{\dummyTypfkAoriginal{\Typpi}}
\newcommand{\dummyTypfkAoriginalpi}[1]{\TypfkAoriginalpi^{#1}}
\newcommand{\TypfkAoriginalpip}{\dummyTypfkAoriginalpi{+}}
\newcommand{\TypfkAoriginalbarpip}{\TypfkAoriginal_{\Typpibar}^{+}}

\newcommand{\TypU}{U}
\newcommand{\TypUchipi}{\TypTHREEdummies{\TypU}{\Typchi}{\Typpi}}

\newcommand{\TyptrK}{K}
\newcommand{\TyptrL}{L}
\newcommand{\TyptrE}{E}
\newcommand{\TyptrF}{F}

\newcommand{\Typtrm}{\varsigma}
\newcommand{\dummyTyptrm}[1]{\Typtrm_{#1}}
\newcommand{\Typtrmone}{\dummyTyptrm{1}}
\newcommand{\Typtrmtwo}{\dummyTyptrm{2}}

\newcommand{\dummyTypU}[1]{\TypU^{#1}}
\newcommand{\TypUo}{\dummyTypU{0}}
\newcommand{\TypUochipi}{\TypTHREEdummies{\TypUo}{\Typchi}{\Typpi}}

\newcommand{\TypUp}{\dummyTypU{+}}
\newcommand{\TypUpchipi}{\TypTHREEdummies{\TypUp}{\Typchi}{\Typpi}}
\newcommand{\TypUm}{\dummyTypU{-}}
\newcommand{\TypUmchipi}{\TypTHREEdummies{\TypUm}{\Typchi}{\Typpi}}

\newcommand{\TypDelta}{\Delta}
\newcommand{\TypS}{S}
\newcommand{\Type}{\varepsilon}
\newcommand{\Typvartheta}{\vartheta}
\newcommand{\Typvarthetachipi}{\TypSecondTHREEdummies{\Typvartheta}{\Typchi}{\Typpi}}

\newcommand{\dummyTypUpflat}[1]{\TypU^{#1}}
\newcommand{\TypUpflat}{\dummyTypUpflat{+,\flat}}
\newcommand{\TypUmflat}{\dummyTypUpflat{-,\flat}}

\newcommand{\TypSh}{{\mathfrak{Sh}}}
\newcommand{\dummyTypSh}[2]{\TypSh^{#1,#2}}
\newcommand{\TypShchipi}{{\dummyTypSh{\Typchi}{\Typpi}}}

\newcommand{\Typkappa}{\kappa}

\newcommand{\TypR}{R}
\newcommand{\TypRchipip}{\TypFOURdummies{\TypR}{\Typchi}{\Typpi}{+}}

\newcommand{\Typvph}{\varphi}

\newcommand{\Typvphpip}{\TypFOURdummies{\Typvph}{\Typchi}{\Typpi}{+}}
\newcommand{\Typomega}{\omega}
\newcommand{\TypprtrZ}{\mathfrak{Z}}
\newcommand{\TypprtrZomchipi}{{\TypTHREEdummies{\TypprtrZ_\Typomega}{\Typchi}{\Typpi}}}

\newcommand{\Typhrho}{{\hat{\rho}}}
\newcommand{\dammyTyphrho}[2]{\Typhrho^{#1,#2}}
\newcommand{\Typhrhochipi}{\dammyTyphrho{\Typchi}{\Typpi}}
\newcommand{\TypfkBchipiomega}{{\TypFOURdummies{{\mathfrak{B}}}{\Typchi}{\Typpi}{\Typomega}}}

\newcommand{\TypfkAoriginalEX}{\dummyTypfkAoriginal{{\mathrm{ex}}}}
\newcommand{\TypfkAoriginalbR}{\dummyTypfkAoriginal{\bR}}
\newcommand{\TypUpsilon}{\Upsilon}
\newcommand{\Typbark}{{\bar k}}

\newcommand{\Typhatdelta}{{\hat{\delta}}}

\newcommand{\Typsgnchipi}{{\mathrm{sgn}}^{\Typchi,\Typpi}}

\newcommand{\TypRchipire}{\TypFOURdummies{\TypR}{\Typchi}{\Typpi}{{\mathrm{real}}}}
\newcommand{\TypRchipinu}{\TypFOURdummies{\TypR}{\Typchi}{\Typpi}{{\mathrm{null}}}}

\newcommand{\Typq}{q}
\newcommand{\dummyTypq}[1]{q_{#1}}
\newcommand{\Typqbeta}{\dummyTypq{\beta}}

\newcommand{\Typc}{c}
\newcommand{\dummyTypc}[1]{\Typc_{#1}}
\newcommand{\Typcbeta}{\dummyTypc{\beta}}
\newcommand{\Typcal}{\dummyTypc{\al}}

\newcommand{\Typomegachhilmb}{{\Typomega^\Typchi_{\lambda,\mu;\beta}}}
\newcommand{\Typhrhochipibeta}{{\Typhrhochipi(\beta)}}
\newcommand{\TypHCbhmpiomega}{{\TypFOURdummies{{\mathfrak{HC}}}{\Typchi}{\Typpi}{\Typomega}}}
\newcommand{\TypmclL}{{\mathcal{L}}}
\newcommand{\TypLam}{\Lambda}

\newcommand{\TypmclLchipiLam}{{\TypSecondFOURdummies{\TypmclL}{\Typchi}{\Typpi}{\TypLam}}}
\newcommand{\TypmclM}{{\mathcal{M}}}
\newcommand{\TypmclMchipiLam}{{\TypSecondFOURdummies{\TypmclM}{\Typchi}{\Typpi}{\TypLam}}}

\newcommand{\TypLamlambda}{\TypLam^{+\lambda}}
\newcommand{\TypLammu}{\TypLam^{+\mu}}
\newcommand{\TypLamnu}{\TypLam^{+\nu}}
\newcommand{\TypLamo}{\TypLam^{+0}}

\newcommand{\Typacsi}{{\acute s}_i}

\newcommand{\TypN}{N}
\newcommand{\TypNchipi}{\TypN^{\Typchi,\Typpi}}

\newcommand{\Typtau}{\tau}
\newcommand{\Typchitau}{\Typtau^\Typchi}
\newcommand{\Typchitauipi}{\Typchitau_i(\Typpi)}

\newcommand{\TypRchipiptaui}{\TypFOURdummies{\TypR}{\Typchi}{\Typchitauipi}{+}}
\newcommand{\TypNchipitaui}{\TypN^{\Typchi,\Typchitauipi}}

\newcommand{\TypT}{T}
\newcommand{\TypTchitauipii}{\TypFOURdummies{\TypT}{\Typchi}{\Typchitauipi}{i}}
\newcommand{\TypUchipitaui}{\TypTHREEdummies{\TypU}{\Typchi}{\Typchitauipi}}

\newcommand{\TyptrEchipi}{\TyptrE^{\Typchi,\Typpi}}
\newcommand{\TyptrFchipi}{\TyptrF^{\Typchi,\Typpi}}

\newcommand{\Typtauichipi}{{\TypFOURdummies{\Typtau}{\Typchi}{\Typpi}{i}}}
\newcommand{\TyptauichipiLam}{\Typtauichipi(\TypLam)}

\newcommand{\TypbBI}{{\mathbb{B}}(\TypfkI)}
\newcommand{\dammyTypchitauTyppi}[2]{\Typtau^\Typchi_{#1,#2}(\Typpi)}
\newcommand{\TypchitauftTyppi}{{\dammyTypchitauTyppi{f}{t}}}
\newcommand{\TypchitauftmTyppi}{{\dammyTypchitauTyppi{f}{t-1}}}
\newcommand{\TypchitaufoTyppi}{{\dammyTypchitauTyppi{f}{0}}}

\newcommand{\TypRchipi}{\TypSecondTHREEdummies{\TypR}{\Typchi}{\Typpi}}

\newcommand{\Typpibar}{{\bar{\Typpi}}}
\newcommand{\Typpialbar}{{\bar{\Typpial}}}
\newcommand{\Typqbar}{{\bar{\Typq}}}
\newcommand{\TypRchibarpi}{\TypSecondTHREEdummies{\TypR}{\Typchi}{\Typpibar}}
\newcommand{\TypRchibarpip}{\TypFOURdummies{\TypR}{\Typchi}{\Typpibar}{+}}
\newcommand{\TypRchibarpipre}{\TypFOURdummies{\TypR}{\Typchi}{\Typpibar}{+,{\mathrm{real}}}}
\newcommand{\TypRchibarpipnu}{\TypFOURdummies{\TypR}{\Typchi}{\Typpibar}{+,{\mathrm{null}}}}
\newcommand{\TypRchibarpire}{\TypFOURdummies{\TypR}{\Typchi}{\Typpibar}{\mathrm{real}}}
\newcommand{\TypRchibarpinu}{\TypFOURdummies{\TypR}{\Typchi}{\Typpibar}{\mathrm{null}}}

\newcommand{\Typs}{s}
\newcommand{\dammyTyps}[1]{\Typs_{#1}}
\newcommand{\Typsbeta}{\dammyTyps{\beta}}
\newcommand{\Typsal}{\dammyTyps{\al}}

\newcommand{\TypW}{W}
\newcommand{\TypWchipire}{\TypFOURdummies{\TypW}{\Typchi}{\Typpi}{{\mathrm{real}}}}

\newcommand{\TypO}{O}
\newcommand{\dammyTypO}[2]{\TypO^{#1,#2}}
\newcommand{\TypOchipi}{\dammyTypO{\Typchi}{\Typpi}}
\newcommand{\TypOchitauipi}{\dammyTypO{\Typchi}{\Typchitauipi}}

\newcommand{\TyphatT}{{\hat{T}}}
\newcommand{\TyphatTchiPRpiPRLami}{\TypFOURdummies{\TyphatT}{\Typchi}{\Typpi^\prime}{\TypLam^\prime,i}}
\newcommand{\TypmcLchiPRpiPRLami}{{\TypSecondFOURdummies{\TypmclL}{\Typchi}{\Typpi^\prime}{\TypLam^\prime}}}
\newcommand{\TypTchiPRpii}{\TypFOURdummies{\TypT}{\Typchi}{\Typpi^\prime}{i}}
\newcommand{\TypUchiPRpi}{\TypTHREEdummies{\TypU}{\Typchi}{\Typpi^\prime}}

\newcommand{\dammyTypchipiLam}[2]{\Typtau^{\Typchi,\Typpi}_{#1,#2}(\TypLam)}
\newcommand{\TypchitauftTyppiLam}{{\dammyTypchipiLam{f}{t}}}
\newcommand{\TypchitauftmTyppiLam}{{\dammyTypchipiLam{f}{t-1}}}
\newcommand{\TypchitaufoTyppiLam}{{\dammyTypchipiLam{f}{0}}}
\newcommand{\TypOchitauipift}{\dammyTypO{\Typchi}{\TypchitauftTyppi}}

\newcommand{\Typhrhochitauipi}{\dammyTyphrho{\Typchi}{\Typchitauipi}}

\newcommand{\TypRchipibeta}{\TypFOURdummies{\TypR}{\Typchi}{\Typpi}{+,\beta}}
\newcommand{\TypRchitauipibeta}{\TypFOURdummies{\TypR}{\Typchi}{\Typchitauipi}{+,\beta}}

\newcommand{\TypP}{P}

\newcommand{\TypPchipipmbeta}{\TypFOURdummies{\TypP}{\Typchi}{\Typpi}{\pm\beta}}
\newcommand{\TypPchipipbeta}{\TypFOURdummies{\TypP}{\Typchi}{\Typpi}{\beta}}
\newcommand{\TypPchipimbeta}{\TypFOURdummies{\TypP}{\Typchi}{\Typpi}{-\beta}}
\newcommand{\TyphatPchipi}{\TypSecondTHREEdummies{{\hat{\TypP}}}{\Typchi}{\Typpi}}

\newcommand{\Typhrhochipifull}{\Typhrhochipi_{{\mathrm{full}}}}

\newcommand{\TypUodag}{\dummyTypU{{0,\dagger}}}
\newcommand{\TypUodagchipi}{\TypUodag(\Typchi,\Typpi)}

\newcommand{\Typhrhodag}{{\hat{\rho}}_\dagger}

\newcommand{\Typsdagbeta}{\Typs^\dagger_\beta}

\newcommand{\Typdots}{{\dot{\Typs}}}
\newcommand{\Typdotsdagbeta}{\Typdots^\dagger_\beta}

\newcommand{\Typepsilon}{\epsilon}
\newcommand{\Typvarpi}{\varpi}

\newcommand{\TypWchipidagre}{\TypFOURdummies{\TypW}{\Typchi}{\Typpi}{\dagger,{\mathrm{real}}}}

\newcommand{\TypdotW}{{\dot{W}}}
\newcommand{\TypdotWchipidagre}{\TypFOURdummies{\TypdotW}{\Typchi}{\Typpi}{\dagger,{\mathrm{real}}}}

\newcommand{\Typxi}{\xi}
\newcommand{\Typxichipidag}{\TypFOURdummies{\Typxi}{\Typchi}{\Typpi}{\dagger}}
\newcommand{\Typdotxi}{{\dot{\Typxi}}}
\newcommand{\Typdotxichipidag}{\TypFOURdummies{\Typdotxi}{\Typchi}{\Typpi}{\dagger}}

\newcommand{\Typomegazero}{{\Typomega_0}}
\newcommand{\TypfkBchipionedagpre}{{\TypFOURdummies{{\mathfrak{B}}}{\Typchi}{\Typpi}{\Typomegazero,\dagger,{\mathrm{pre}}}}}
\newcommand{\TypfkBchipiomegazero}{{\TypFOURdummies{{\mathfrak{B}}}{\Typchi}{\Typpi}{\Typomegazero}}}

\newcommand{\TypmclE}{{\mathfrak{E}}}
\newcommand{\dammyTypmclE}[1]{\TypmclE^{#1}}
\newcommand{\TypmclELam}{\dammyTypmclE{\TypLam}}
\newcommand{\TypeLam}{e^\TypLam}

\newcommand{\TypRchipipre}{\TypFOURdummies{\TypR}{\Typchi}{\Typpi}{+,{\mathrm{real}}}}
\newcommand{\TypRchipipnu}{\TypFOURdummies{\TypR}{\Typchi}{\Typpi}{+,{\mathrm{null}}}}

\newcommand{\TypmBiclELam}{{\mathrm{Aut}}_\bZ(\TypmclELam)}%{{\mathrm{Bi}}(\TypmclELam)}
\newcommand{\TypsigLam}{\sigma^\TypLam}
\newcommand{\TypfkWchipireLam}{{\mathfrak{W}}^{\Typchi,\Typpi,\TypLam}_{{\mathrm{real}}}}
\newcommand{\TypdotsigLam}{{\dot \sigma}^\TypLam}
\newcommand{\TypdotfkWchipireLam}{{\dot{{\mathfrak{W}}}}^{\Typchi,\Typpi,\TypLam}_{{\mathrm{real}}}}

\newcommand{\TypBifkApi}{{\mathrm{Bi}}(\TypfkAoriginalpi)}
\newcommand{\TypsLam}{s^\TypLam}
\newcommand{\TypdotsLam}{{\dot s}^\TypLam}
\newcommand{\TypWchipireLam}{W^{\Typchi,\Typpi,\TypLam}_{{\mathrm{real}}}}
\newcommand{\TypdotWchipireLam}{{\dot{W}}^{\Typchi,\Typpi,\TypLam}_{{\mathrm{real}}}}

\newcommand{\TypxichipiLam}{\TypFOURdummies{\Typxi}{\Typchi}{\Typpi}{\TypLam}}
\newcommand{\TypdotxichipiLam}{\TypFOURdummies{\Typdotxi}{\Typchi}{\Typpi}{\TypLam}}

\newcommand{\TypetachipiLamp}{\eta^{\Typchi,\Typpi}_{\TypLam,\prime}}
\newcommand{\TypdotetachipiLamp}{{\dot{\eta}}^{\Typchi,\Typpi}_{\TypLam,\prime}}

\newcommand{\TypetachipiLam}{\eta^{\Typchi,\Typpi}_\TypLam}
\newcommand{\TypdotetachipiLam}{{\dot{\eta}}^{\Typchi,\Typpi}_\TypLam}

\newcommand{\TyphatmclELam}{{\hat{\TypmclE}}^\TypLam}
\newcommand{\TypcheckmclELam}{{\check{\TypmclE}}^\TypLam}

\newcommand{\TypiotaLam}{\iota^\TypLam}

\newcommand{\TypmclEomegazero}{\dammyTypmclE{\Typomegazero}}
\newcommand{\Typeomegazero}{e^{\Typomegazero}}

\newcommand{\acutechipi}{{{\acute{\Typchi}},{\acute{\Typpi}}}}

\begin{document}

\begin{spacing}{2}
\begin{center}
{\huge{
Typical Irreducible Characters of}} \\
{\huge{Generalized Quantum Groups}}
\end{center}
\end{spacing}

%\vspace{1mm}

%\begin{center}
%{\it{Dedicated to Professor Nicol{\'{a}}s Andruskiewitsch \\
%on the occasion of his 60th birthday}}
%\end{center}

%\vspace{2mm}

\begin{center}
%by \\ 
Hiroyuki Yamane\footnote{This work was supported by
JSPS Grant-in-Aid for Scientific Research (C) 19K03420.}
\end{center}

\begin{abstract}
We introduce the
definition of the typical irreducible modules
of the generalized quantum groups,
and prove the Weyl-Kac-type formulas of 
their characters.
As a by-product,
we obtain the Weyl-Kac-type  character formulas of 
the typical irreducible modules
of the quantum superalgebras
associated with the basic classical Lie superalgebras,
which is explained in Introduction.
\end{abstract}

\section{Introduction}
In 1977, V.~Kac~\cite{Kac77b}\cite{Kac77c}  gave the Weyl-type character formula of 
the typical finite-dimensional irreducible modules
of the basic classical Lie superalgebras ${\mathfrak{k}}$.
See also \cite[Subsection~8.6 and Theorem~14.4.2]{Musson12}.
The highest weights of the non-typical (atypical) finite-dimensional irreducible modules
of ${\mathfrak{k}}$
stay in a Zariskii closed subset of the linear space
spanned by those of all the finite-dimensional irreducible modules of ${\mathfrak{k}}$.
In this paper, we give its counterpart
for the generalized quantum groups
over any field $\bK$ of characteristic zero.
As for the N.~Geer's related result \cite{Geer07}, see Remark~\ref{remark:Geer}.

Let $\bC$ be the field of complex numbers.
From now on until the end of Introduction, associative algebras and Lie superalgerbas 
are those over $\bC$.
Let $\Typl\in\bN$ and $I:=\{1,\ldots,\Typl\}$.
Let $q_{ij}\in\bCt:=\bC\setminus\{0\}$
$(i, j\in I)$, and $\TypIntQ:=(q_{ij})_{i,j\in I}$ (the $\Typl\times\Typl$-matrix).
In the same way as that for the Lusztig's definition of the quantum groups
in \cite[Subsection~3.1.1]{b-Lusztig93},
we define the generalized quantum groups
$U_\TypIntQ$. As an associative $\bC$-algebra,
$U_\TypIntQ$ has generators
$K_i^{\pm 1}$, $L_i^{\pm 1}$, $E_i$, $F_i$ $(i\in\TypfkI)$ satisfying the equations
\begin{equation*}
\begin{array}{l}
K_iK_j=K_jK_i,\,L_iL_j=L_jL_i,\,K_iL_j=L_jK_i,\\
K_iK_i^{-1}=K_i^{-1}K_i=L_iL_i^{-1}=L_i^{-1}L_i=1, \\
K_iE_i=q_{ij}E_jK_i,\,q_{ij}K_iF_i=F_jK_i,\,q_{ji}L_iE_i=E_jL_i,\,L_iF_i=q_{ji}F_jL_i,\\
E_iF_j-F_jE_i=\delta_{ij}(-K_i+L_i),
\end{array}
\end{equation*} $(i,j\in I)$.
Let $U_\TypIntQ^0$ (resp. $U_\TypIntQ^+$,
resp. $U_\TypIntQ^-$) be the $\bC$-subalgebras of $U_\TypIntQ$
generated by $\{K_i^{\pm 1},L_i^{\pm 1}|i\in I\}$
(resp.  $\{E_i^{\pm 1}|i\in I\}$, resp.  $\{F_i^{\pm 1}|i\in I\}$). 
Then we have the $\bC$-linear isomorphism $g:U_\TypIntQ^-\otimes 
U_\TypIntQ^0\otimes U_\TypIntQ^+\to U_\TypIntQ$ with $g(Y\otimes Z\otimes X):=YZX$,
we identify $U_\TypIntQ^0$
with the Laurent polynomial $\bC$-algebra $\bC[K_i^{\pm 1},L_i^{\pm 1}|i\in I]$,
and we idnetify $U_\TypIntQ^+$ and $U_\TypIntQ^-$
with the Nichols algebra of diagonal-type with the braidings defined from
$\TypIntQ$ and $\TypIntQ^{\mathrm{op}}:=(q_{ji})_{i,j\in I}$
respectively. 
Let $\TypIntKhR^\TypIntQ_+$ be the Kharchenko's 
positive root system of $U_\TypIntQ^+$ (and $U_\TypIntQ^-$).
In this paper, we introduce an appropriate definition of 
the typical finite-dimensional irreducible modules of $U_\TypIntQ$ in the case where $|\TypIntKhR^\TypIntQ_+|
<\infty$,
and prove the Weyl-Kac type forumula of their characters.
For such $U_\TypIntQ$, 
we have already achieved 
the (axiomatic) Weyl groupoids \cite{HY08}, \cite{CH15},
the Bruhat order of the Weyl groupoids \cite{AY18},
the Shapovalov determinants \cite{HY10},
the defining relations \cite{A15} (see also \cite{AAY11}, \cite{JMY18}, \cite{Y94}, \cite{Y99}, \cite{Y01}),
the universal $R$-matrices \cite{AY15},
the classification of the finite-dimensional irreducible modules \cite{AYY15} 
and the Harish-Chandra-type theorem for the (skew) center \cite{BY18}
(see also \cite{BY15}, \cite{BY19}).

Let ${\mathfrak{a}}$ be a linear space over $\bC$.
Let ${\mathfrak{a}}(0)$ and ${\mathfrak{a}}(1)$
be linear subspaces of ${\mathfrak{a}}$
such that ${\mathfrak{a}}={\mathfrak{a}}(0)\oplus{\mathfrak{a}}(1)$.
For $t\in\bZ$, 
define ${\mathfrak{a}}(t):={\mathfrak{a}}(t\pm 2)$.
Let $[\,,\,]:{\mathfrak{a}}\times{\mathfrak{a}}\to{\mathfrak{a}}$
be a $\bC$-bilinear map with $[{\mathfrak{a}}(t),{\mathfrak{a}}(t^\prime)]
\subset{\mathfrak{a}}(t+t^\prime)$.
We say that ${\mathfrak{a}}=({\mathfrak{a}},[\,,\,])$
is {\it{a Lie superalgebra over $\bC$}} if
\begin{equation*}
\begin{array}{l}
[X,[Y,Z]]=[[X,Y],Z]+(-1)^{st}[Y,[X,Z]],\quad [Y,X]=- (-1)^{st}[X,Y]
\\
 (X\in{\mathfrak{a}}(s),\,
Y\in{\mathfrak{a}}(t),\,Z\in{\mathfrak{a}}).
\end{array}
\end{equation*}

Let ${\mathfrak{b}}$ be one of the Lie superalgebras over $\bC$
which are 
$A(m-1,n-1)={\mathfrak{sl}}(m|n)$ $(\Typl=m+n-1\geq 3,1\leq m<n)$,
${\mathfrak{sl}}(n|n)$ $(\Typl=2n-1\geq 3)$,
$B(m,n)={\mathfrak{osp}}(2m+1|2n)$ $(\Typl=m+n\geq 1, m\geq 0,  n\geq 1)$,
$C(n)={\mathfrak{osp}}(2|2n-2)$ $(\Typl=n\geq 3)$,
$D(m,n)={\mathfrak{osp}}(2m|2n)$ $(\Typl=m+n\geq 3, m\geq 2,  n\geq 1)$,
$F(4)$ $(\Typl=4)$,
$G(3)$ $(\Typl=3)$,
$D(2,1;a)$ $(\Typl=3, a\in\bC\setminus\{0,-1\})$,
where $m$, $n\in\bZgeqo$.
We call these Lie superalgebras
{\it{the finite-dimensional Lie superalgebras of type $A$-$G$
of rank $=\Typl$}}.
The family of {\it{the basic classical Lie superalgebras of rank $=\Typl$}}
is composed of 
these Lie superalgebras except for ${\mathfrak{sl}}(n|n)$ $(\Typl+1\in 2\bZ)$
and $A(n-1,n-1)={\mathfrak{sl}}(n|n)/\bC {\mathcal I}$ $(\Typl=2n-1\geq 3)$,
where ${\mathcal I}$ is the $2n\times 2n$ identity matrix.
Each of this family is a simple Lie superalgebra.

The quantum superalgebra $U_q({\mathfrak{b}}$) of
the above ${\mathfrak{b}}$ can be obtained from a generalized quantum group
by slight modification (see Remark~\ref{remark:relGQGandSUP}).
As Theorem~\ref{theorem:IntroMain} below,
we state a formula of the characters of the typical finite-dimensional irreducible modules of  $U_q({\mathfrak{b}})$.

Let ${\mathcal{E}}$ be a $\Typl$-dimensional linear space over $\bC$
with a $\bC$-basis $\Pi:=\{\al_i|i\in I\}$.
Let $(\,,\,):{\mathcal{E}}\times{\mathcal{E}}\to\bC$ be a symmetric bilinear map.
Let $p:\bZ\Pi\to\bZ$ be a $\bZ$-module homomorphism 
such that $p(\Pi)\subset\{0,1\}$.
We define the Lie superalgebra
${\mathfrak{g}}={\mathfrak{g}}((\,,\,),\Pi,p)={\mathfrak{g}}(0)\oplus{\mathfrak{g}}(1)$ over $\bC$ by the following
axioms $({\mathfrak{g}}1)$-$({\mathfrak{g}}3)$. \newline\newline
$({\mathfrak{g}}1)$ There exists a direct sum 
${\mathfrak{g}}=\oplus_{\lambda\in\bZ\Pi}{\mathfrak{g}}_\lambda$
as a $\bC$-linear space such that $[{\mathfrak{g}}_\lambda,{\mathfrak{g}}_\mu]
\subset{\mathfrak{g}}_{\lambda+\mu}$ 
$(\lambda,\mu\in\bZ\Pi)$ and ${\mathfrak{g}}(t)=\oplus_{p(\lambda)-t\in 2{\mathbb{Z}}}{\mathfrak{g}}_\lambda$
$(t\in\{0,1\})$.
\newline 
$({\mathfrak{g}}2)$ There exist $e_i\in{\mathfrak{g}}_{\al_i}$, $f_i\in{\mathfrak{g}}_{-\al_i}$
$(i\in\TypfkI)$
and
a $\bC$-basis $\{h_i|i\in I\}$ of ${\mathfrak{g}}_0$
such that 
\begin{equation*}
[h_i,h_j]=0,\,[h_i,e_j]=(\al_i,\al_j)e_j,\,[h_i,f_j]=-(\al_i,\al_j)f_j,\,[e_i,f_j]=\delta_{ij}h_j\quad (i,j\in I).
\end{equation*}
(Hence $\dim{\mathfrak{g}}_0=|I|$, $\dim{\mathfrak{g}}_{\pm\al_i}=1$, 
and $\{\lambda\in\bZ\Pi|\dim{\mathfrak{g}}_\lambda\geq 1\}\subset
\bZgeqo\Pi\cup\bZleqo\Pi$.)
\newline 
$({\mathfrak{g}}3)$ 
\begin{equation*}
\begin{array}{l}
\mbox{$\forall\lambda\in\bZgeqo\Pi\setminus\{0\}$ (resp. $\bZleqo\Pi\setminus\{0\}$),} \\
\mbox{$\{X\in{\mathfrak{g}}_\lambda|\,\forall i\in I,[f_i,X]=0\,\mbox{(resp. $[e_i,X]=0$)}\}=\{0\}$.}
\end{array}
\end{equation*}
\newline
Let $\TypIntRp:=\{\beta\in\bZgeqo\Pi\setminus\{0\}|\dim{\mathfrak{g}}_\beta\geq 1\}$,
$\TypIntRpnull:=\{\beta\in\TypIntRp|(\beta,\beta)=0\}$
and 
$\TypIntRpreal:=\{\beta\in\TypIntRp|(\beta,\beta)\ne 0, \beta\notin 2\TypIntRp\}$.
If $\dim{\mathfrak{g}}<\infty$, 
$\TypIntRp=\TypIntRpnull\cup\TypIntRpreal
\cup\{2\beta|\beta\in\TypIntRpreal, p(\beta)\in 2\bZ\}$.
Let $\TypIntR:=\TypIntRp\cup(-\TypIntRp)$.
Then ${\mathfrak{g}}={\mathfrak{h}}\oplus(\oplus_{\beta\in\TypIntR}{\mathfrak{g}}_\beta)$,
where let ${\mathfrak{h}}:={\mathfrak{g}}_0$.
Let $(\TypIntRpnull)_t:=\{\gamma\in\TypIntRpnull|(-1)^{p(\gamma)}=(-1)^t\}$
$(t\in\{0,1\})$.
Let $\bCt:=\bC\setminus\{0\}$
and $q\in\bCt\setminus\{\pm 1\}$.
Let $a\in\bC$ be such that $q=\exp(a)$.
For $b\in\bC$, let $q^b:=\exp(ab)$ as usual.
Let 
\begin{equation*}
{\acute{q}}_{ij}:=(-1)^{p(\al_i)p(\al_j)}q^{(\al_i,\al_j)}(\in\bCt)\quad(i,j\in I).
\end{equation*}  
Let $\TypIntacuteQ:=({\acute{q}}_{ij})_{i,j\in I}$.
Let $\bCtinf:=\{z\in\bCt|\forall n\in\bN,z^n\ne 1\}$.
Assume that $\prod_{t=1}^k{\acute{q}}_{i_tj_t}\in\bCtinf$ 
for
$k\in\bN$ and
$i_t$, $j_t\in I$ $(1\leq t\leq k)$
with
$\prod_{t=1}^k(\al_{i_t},\al_{j_t})\ne 0$. 
Then if $\dim{\mathfrak{g}}<\infty$, 
$\TypIntKhR^\TypIntacuteQ_+$ can be identified with
$\TypIntRpnull\cup\TypIntRpreal$
(see Remark~\ref{remark:relofInt}).

Assume $\Typl=1$. Then $\dim{\mathfrak{g}}<\infty$.
In fact, ${\mathfrak{g}}$ is isomorphic to
${\mathfrak{sl}}_2(\bC)$
(resp. ${\mathfrak{osp}}(1|2)$, resp. ${\mathfrak{sl}}(1|1)$,
resp. ${\mathfrak{h}}_1$)
if $p(\al_1)=0$ (resp. $1$, resp. $1$, resp. $0$)
and $(\al_1,\al_1)\ne 0$ 
(resp. $\ne 0$, resp. $=0$, resp. $=0$),
where ${\mathfrak{h}}_1$ is the Heisenberg Lie algebra.

Assume that $\Typl\geq 2$, $\dim{\mathfrak{g}}<\infty$ and 
there exists $(i,j)\in I^\prime\times(I\setminus I^\prime)$
with $(\al_i,\al_j)\ne 0$
for every non-empty proper subset $I^\prime$ of $I$.
We have $|\TypIntRp|<\infty$ and $(\TypIntRpnull)_0=\emptyset$ for $U_{{\acute{Q}}}$.
The Weyl groupoids of ${\mathfrak{g}}$ and $U_{{\acute{Q}}}$ are isomorphic
(see \cite[Lamma~5.1~(2)]{AYY15}).
From this fact, 
using the classification theorem \cite{Hec09} of the (finite-type) Nichols algebras of 
diagonal-type, we obtain a new proof that
${\mathfrak{g}}$ is isomorphic to one of
the finite dimensional simple Lie algebras of rank~$=\Typl$
and the finite-dimensional Lie superalgebras of type $A$-$G$
of rank~$=\Typl$, which has originally been proved by \cite{Kac77a}.

Assume $\dim{\mathfrak{g}}<\infty$ until the end of Introduction.
Let $U_q({\mathfrak{g}})$ be the quantum superalgebra
of ${\mathfrak{g}}$. As an associative $\bC$-algebra,
$U_q({\mathfrak{g}})$ is defined with the
generators $k_i^{\pm 1}$, $e_i$, $f_i$ $\in U_q({\mathfrak{g}})$ $(i\in I)$
and the relations
composed of those of (QS1-3) below and the same ones as those of 
\cite[(QS4)(1-17), (QS5)(1-17) of Proposition~6.7.1]{Y99}
(see also \cite[Proposition~10.4.1]{Y94}).
\newline\newline
(QS1) \hspace{2cm} $k_ik_i^{-1}=k_i^{-1}k_i=1,\,\,k_ik_j=k_ik_j,$ \newline
(QS2) \hspace{2cm} $k_ie_jk_i^{-1}=q^{(\al_i,\al_j)}e_j,\,\,k_if_jk_i^{-1}=q^{-(\al_i,\al_j)}f_j$, \newline
(QS3) \hspace{2cm} $e_if_j-(-1)^{p(\al_i)p(\al_j)}f_je_i=\delta_{ij}{\frac {k_i-k_i^{-1}} {q-q^{-1}}}$.
\newline\newline
We have the $\bC$-linear subspaces $U_q({\mathfrak{g}})_\lambda$
$(\lambda\in\bZ)$ of $U_q({\mathfrak{g}})$
such that $k_i^{\pm}\in U_q({\mathfrak{g}})_0$,
$e_i\in U_q({\mathfrak{g}})_{\al_i}$,
$f_i\in U_q({\mathfrak{g}})_{-\al_i}$ $(i\in I)$
and $U_q({\mathfrak{g}})=\oplus_{\lambda\in\bZ\Pi}U_q({\mathfrak{g}})_\lambda$.

Let ${\mathrm{Gr}}_\Pi$ be the set of group homomorphisms
from $\bZ\Pi$ to $\bCt$.
Define ${\hat{\rho}}_\Pi\in{\mathrm{Gr}}_\Pi$ by
${\hat{\rho}}_\Pi(\al_i):={\acute{q}}_{ii}$ $(i\in I)$.
It follows from Proposition~\ref{proposition:mqbeta} that
\begin{equation*}
\forall\beta\in\TypIntRpreal,\,\,\exists r^{{\hat{\rho}}}(\beta)\in\bZ,\,\,{\rm{s.t.}}\,\,
{\hat{\rho}}_\Pi(\beta)=((-1)^{p(\beta)}q^{(\beta,\beta)})^{r^{{\hat{\rho}}}(\beta)}.
\end{equation*}

Let ${\mathrm{Bi}}_\Pi$ be the group of 
bijections from $\bZ\Pi$ to $\bZ\Pi$,
where the multiplication is define by $st:=s\circ t$ $(s,t\in{\mathrm{Bi}}_\Pi)$.

Let $\Omega\in{\mathrm{Gr}}_\Pi$.
Let $\TypIntmcMOm$ (resp. $\TypIntmcLOm$)
be the Verma module (resp. the highest weight irreducible module) of $U_q({\mathfrak{g}})$
with the highest weight vector ${\tilde{v}}_\Omega$
(resp. $v_\Omega$)
satisfying  $k_i^{\pm 1}{\tilde{v}}_\Omega=\Omega(\pm\al_i){\tilde{v}}_\Omega$
and $e_i{\tilde{v}}_\Omega=0$
(resp.  $k_i^{\pm 1}v_\Omega=\Omega(\pm\al_i)v_\Omega$
and $e_iv_\Omega=0$).
Let $\TypIntmcMOm_\lambda:=U_q({\mathfrak{g}})_\lambda{\tilde{v}}_\Omega$
and $\TypIntmcLOm_\lambda:=U_q({\mathfrak{g}})_\lambda v_\Omega$ $(\lambda\in\bZ\Pi)$.
Then $\TypIntmcMOm=\oplus_{\lambda\in\bZleqo\Pi}\TypIntmcMOm_\lambda$,
$\TypIntmcLOm=\oplus_{\lambda\in\bZleqo\Pi}\TypIntmcLOm_\lambda$
and $\dim\TypIntmcMOm_\mu=\dim\TypIntmcLOm_\mu=0$
$(\mu\in\bZ\Pi\setminus\bZleqo\Pi)$. 
For $\lambda\in\bZ\Pi$,
$\dim\TypIntmcMOm_\lambda$ equals the cardinality
of the set of maps $f:\TypIntRpreal\cup\TypIntRpnull\to\bZgeqo$
with $f((\TypIntRpnull)_1)\subset\{0,1\}$
and $\sum_{\beta\in\TypIntRpreal\cup\TypIntRpnull}f(\beta)\beta=\lambda$.
For $\nu\in\bZ\Pi$,
define $\Omega^{+\nu}\in{\mathrm{Gr}}_\Pi$
by $\Omega^{+\nu}(\lambda):=(-1)^{p(\nu)p(\lambda)}q^{(\nu,\lambda)}\Omega(\lambda)$ $(\lambda\in\bZ\Pi)$.

Assume $\dim\TypIntmcLOm<\infty$.
It follows that
\begin{equation*}
\forall\beta\in\TypIntRpreal,\,\,\exists n^\Omega(\beta)\in\bZgeqo,\,\,{\rm{s.t.}}\,\,
\Omega(2\beta)=((-1)^{p(\beta)}q^{(\beta,\beta)})^{n^\Omega(\beta)}.
\end{equation*}
For $\beta\in\TypIntRpreal$,
define $s_\beta$, ${\dot{s}}^\Omega_\beta\in{\mathrm{Bi}}_\Pi$ by 
\begin{equation*}
s_\beta(\lambda):=\lambda-{\frac {2(\beta,\lambda)} {(\beta,\beta)}}\beta,\quad
{\dot{s}}^\Omega_\beta(\lambda):=s_\beta(\lambda)
-(r^{{\hat{\rho}}}(\beta)+n^\Omega(\beta))\beta\quad(\lambda\in\bZ\Pi).
\end{equation*}
Let $W^\Pi$ (resp. ${\dot{W}}^\Omega$) be the subgroup of 
${\mathrm{GL}}({\mathcal{E}})$ (resp. ${\mathrm{Bi}}_\Pi$) 
generated by $\{s_\beta|\beta\in\TypIntRpreal\}$
(resp. $\{{\dot{s}}^\Omega_\beta|\beta\in\TypIntRpreal\}$).
Then we have the group isomorphism $h: W^\Pi\to{\dot{W}}^\Omega$
defined by $h(s_\beta):={\dot{s}}^\Omega_\beta$ $(\beta\in\TypIntRpreal)$.
There exists a group homomorphism
${\mathrm{sgn}}^\Omega:{\dot{W}}^\Omega\to\{\pm 1\}$
such that ${\mathrm{sgn}}^\Omega({\dot{s}}^\Omega_\beta)=-1$.

By our main theorem Theorem~\ref{theorem:ScMain} below, we have:
\begin{theorem}\label{theorem:IntroMain}
Assume that 
\begin{equation}\label{eqn:AssScMain}
\mbox{$(\TypIntRpnull)_0=\emptyset$,
$\dim\TypIntmcLOm<\infty$
and $\prod_{\gamma\in\TypIntRpnull}(1-(-1)^{p(\gamma)}{\hat{\rho}}_\Pi(\gamma)\Omega(2\gamma))\ne 0$.}
\end{equation}
Then we have
\begin{equation}\label{eqn:FormulaScMain}
\dim\TypIntmcLOm_\lambda=\sum_{{\dot{w}}\in{\dot{W}}^\Omega}
{\mathrm{sgn}}^\Omega({\dot{w}})\cdot\dim{\mathcal{M}}(\Omega^{+{\dot{w}}(0)})_{\lambda-{\dot{w}}(0)}
\quad(\lambda\in\bZ\Pi).
\end{equation}
\end{theorem}

We call $\TypIntmcLOm$ of \eqref{eqn:AssScMain} {\it{typical}}.
The formula \eqref{eqn:FormulaScMain} for $\TypIntRpnull=\emptyset$
has been obtained by \cite[Theorem~4.15~(3)]{ASR10}.

\begin{remark}\label{remark:relGQGandSUP}
Since $\TypIntacuteQ$ is a $\Typl\times\Typl$ symmetric matrix,
for every $i\in\TypfkI$,  $K_iL_i-1$ belongs to the center
$Z({\bar U}_\TypIntacuteQ)$ of ${\bar U}_\TypIntacuteQ$.
Consider the quotient $\bC$-algebra
${\bar U}_\TypIntacuteQ:=U_\TypIntacuteQ/(\sum_{i\in I}(K_iL_i-1)U_\TypIntacuteQ)$.
Let $d:U_\TypIntacuteQ:={\bar U}_\TypIntacuteQ$ be the canonical map.
Let ${\bar K}_i:=d(K_i)$, ${\bar E}_i:=d(E_i)$, ${\bar F}_i:=d(F_i)$.
Let ${\bar U}_\TypIntacuteQ^\sigma={\bar U}_\TypIntacuteQ\oplus {\bar U}_\TypIntacuteQ\sigma$ 
(resp. $U_q({\mathfrak{g}})^\sigma=U_q({\mathfrak{g}})\oplus U_q({\mathfrak{g}})\sigma$) be the associative 
$\bC$-algebra obtained from ${\bar U}_\TypIntacuteQ$ (resp. $U_q({\mathfrak{g}})$)
by adding the element $\sigma$ with $\sigma^2=1$,
$\sigma {\bar K}_i^{\pm 1}\sigma={\bar K}_i^{\pm 1}$,
$\sigma {\bar E}_i\sigma=(-1)^{p(\al_i)}{\bar E}_i$,
$\sigma {\bar F}_i\sigma=(-1)^{p(\al_i)}{\bar F}_i$
(resp. $\sigma k_i^{\pm 1}\sigma=k_i^{\pm 1}$
$\sigma e_i\sigma=(-1)^{p(\al_i)}e_i$,
$\sigma f_i\sigma=(-1)^{p(\al_i)}f_i$).
Then we have the $\bC$-algebra isomorphism
$t:{\bar U}_\TypIntacuteQ^\sigma\to U_q({\mathfrak{g}})^\sigma$ defined by
$t({\bar K}_i^{\pm 1}):=k_i^{\pm 1}\sigma^{p(\al_i)}$,
$t({\bar E}_i):=e_i$,
$t({\bar F}_i):=-(q-q^{-1})f_i\sigma^{p(\al_i)}$, $t(\sigma):=\sigma$.
Moreover 
$\TypIntKhR^\TypIntacuteQ_+$
can be identified with $\TypIntRp\setminus 2\TypIntRp$.
See also Remark~\ref{remark:relofInt}.
\end{remark}

\section{Generalized guantum groups}
\subsection{Preliminaries}

For $x$, $y\in\bR\cup\{\pm\infty\}$,
let $\TypfkJ_{x,y}:=\{z\in\bZ|x\leq z\leq y\}$.
Let $\bK$ be a field. Let $\bKt:=\bK\setminus\{0\}$.
For $n\in\bZgeqo$ and $x\in\bK$, let $(n)_x:=\sum_{r=1}^n x^{r-1}$,
and $(n)_x!:=\prod_{r=1}^n(r)_x$.
%For $n\in\bZgeqo$, $m\in\TypfkJ_{0,n}$ and $x\in\bK$,
%define ${n\choose m}_x\in\bK$ by
%${n\choose 0}_x:={n\choose n}_x:=1$,
%and ${n\choose m}_x:={n-1\choose m}_x+x^{n-m}{n-1\choose m-1}_x
%=x^m{n-1\choose m}_x+{n-1\choose m-1}_x$
%(if $m\in\TypfkJ_{1,n-1}$).
%If $(m)_x!(n-m)_x!\ne 0$, then ${n\choose m}_x=
%{\frac {(n)_x!} {(m)_x!(n-m)_x!}}$.
%For $x$, $y$, $z\in\bK$, and $n\in\bN$, we have
%$\prod_{t=0}^{n-1}(y+x^tz)=\sum_{m=0}^n
%x^{{\frac {m(m-1)} 2}}{n\choose m}_xy^{n-m}z^m$.
For $n\in\bZgeqo$ and $x$, $y\in\bKt$,
let $(n;x,y):=1-x^{n-1}y$ and $(n;x,y)!:=\prod_{m=1}^n(m;x,y)$.

For $x\in\bKt$, define $\Typkpch(x)\in\bZgeqo\setminus\{1\}$ by
\begin{equation}\label{eqn:chd}
\Typkpch(x):= 
\left\{\begin{array}{l}\min\{\,r^\prime\in\TypfkJ_{2,\infty}\,|\,(r^\prime)_x!=0\,\} 
\,\,\mbox{if $(r^{\prime\prime})_x!=0$ for some $r^{\prime\prime}\in\TypfkJ_{2,\infty}$}, \\
0 \quad\mbox{otherwise.}
\end{array}\right.
\end{equation} 

Let $\TypfkAoriginal$ be a non-zero finite rank free $\bZ$-module.
Let $\Typlp:={\mathrm{Rank}}_\bZ(\TypfkAoriginal)(\in\bN)$.
Let $\Typl\in\TypfkJ_{1,\Typlp}$ and $\TypfkI:=\TypfkJ_{1,\Typl}$.
Let $\{\Typpial_i|i\in\TypfkI\}\cup\{\Typepsilon_r|r\in\TypfkJ_{\Typl+1,\Typlp}\}$ be a $\bZ$-base of $\TypfkAoriginal$.
Define an injection $\Typpi:\TypfkI\to\TypfkAoriginal$ 
by $\Typpial_i:=\Typpi(i)$ ($i\in\TypfkI$).
Let $\TypfkAoriginalpi=\oplus_{i\in\TypfkI}\bZ\Typpial_i$
and $\TypfkAoriginalpip:=\oplus_{i\in\TypfkI}\bZgeqo\Typpial_i$.
Let $\TypfkAoriginalEX:=\oplus_{r=1}^{\Typlp-\Typl}\bZ\Typepsilon_r$.
Then $\TypfkAoriginal=\TypfkAoriginalpi\oplus\TypfkAoriginalEX$.
Note that $\TypfkAoriginalEX=\{0\}$ if $\Typlp=\Typl$.

Let $\Typchi:\TypfkAoriginal\times\TypfkAoriginal\to\bKt$ be a map such that
\begin{equation}\label{eqn:bich}
\forall \lambda, \forall \mu, \forall \nu\in\TypfkAoriginal,\,\Typchi(\lambda+\mu,\nu)=\Typchi(\lambda,\nu)\Typchi(\mu,\nu),\,
\Typchi(\lambda,\mu+\nu)=\Typchi(\lambda,\mu)\Typchi(\lambda,\nu).
\end{equation}

There exists a unique associative $\bK$-algebra (with $1$)
$\TypU=\TypUchipi$ satisfying the following conditions $(\TypU 1)$-$(\TypU 6)$.
\newline\newline
$(\TypU 1)$ As a $\bK$-algebra, $\TypU$ is generated by the elements:
\begin{equation}\label{eqn:gene}
\TyptrK_\lambda,\,\TyptrL_\lambda\,\,(\lambda\in\TypfkAoriginal),\quad
\TyptrE_i, \TyptrF_i\,\,(i\in\TypfkI).
\end{equation} \newline
$(\TypU 2)$ The elements of \eqref{eqn:gene} satisfy the following relations.
\begin{equation}\label{eqn:relone}
\begin{array}{l}
\TyptrK_0=\TyptrL_0=1,\,
\TyptrK_\lambda \TyptrK_\mu=\TyptrK_{\lambda+\mu},\,
\TyptrL_\lambda \TyptrL_\mu=\TyptrL_{\lambda+\mu},\,
\TyptrK_\lambda \TyptrL_\mu=\TyptrL_\mu \TyptrK_\lambda, \\
\TyptrK_\lambda \TyptrE_i =\Typchi(\lambda,\Typpial_i)\TyptrE_i \TyptrK_\lambda,\,
\TyptrL_\lambda \TyptrE_i  =\Typchi(-\Typpial_i,\lambda)\TyptrE_i \TyptrL_\lambda,\\
\TyptrK_\lambda \TyptrF_i  =\Typchi(\lambda,-\Typpial_i)\TyptrF_i \TyptrK_\lambda,\,
\TyptrL_\lambda \TyptrF_i  =\Typchi(\Typpial_i,\lambda)\TyptrF_i\TyptrL_\lambda,\\
\mbox{$\TyptrE_i\TyptrF_j-\TyptrF_j\TyptrE_i=\delta_{ij}(-\TyptrK_{\Typpial_i}+\TyptrL_{\Typpial_i})$}.
\end{array}
\end{equation}\newline
$(\TypU 3)$ Define the map $\Typtrmone:\TypfkAoriginal\times\TypfkAoriginal\to\TypU$
by $\Typtrmone(\lambda,\mu):=\TyptrK_\lambda\TyptrL_\mu$.
Define the $\bK$-subalgebra $\TypUo=\TypUochipi$ of $\TypU$ by 
$\TypUo:=\TyprmSpan_\bK(\Typtrmone(\TypfkAoriginal\times\TypfkAoriginal))$.
Then  $\Typtrmone$ is injective,
and $\Typtrmone(\TypfkAoriginal\times\TypfkAoriginal)$ is a $\bK$-basis of  $\TypUo$.
\newline\newline
$(\TypU 4)$ There exist the $\bK$-subspaces $\TypU_\lambda=\TypUchipi_\lambda$ of $\TypU$
$(\lambda\in\TypfkAoriginalpi)$ satisfying the following conditions $(\TypU 4{\mbox{-}}1)$-$(\TypU 4{\mbox{-}}3)$.
\newline\newline
$(\TypU 4{\mbox{-}}1)$ We have $\TypUo\subset\TypU_0$ and $\TyptrE_i\in\TypU_{\Typpial_i}$,
$\TyptrF_i\in\TypU_{-\Typpial_i}$ ($i\in\TypfkI$). \newline
$(\TypU 4{\mbox{-}}2)$ We have $\TypU_\lambda\TypU_\mu\subset
\TypU_{\lambda+\mu}$ ($\lambda$, $\mu\in\TypfkAoriginalpi$).  \newline
$(\TypU 4{\mbox{-}}3)$ We have 
$\TypU=\oplus_{\lambda\in\TypfkAoriginalpi}\TypU_\lambda$
as a $\bK$-linear spaces.
\newline\newline
$(\TypU 5)$ Let $\TypUp=\TypUpchipi$ 
(resp. $\TypUm=\TypUmchipi$) be the $\bK$-subalgebra (with $1$) of $\TypU$
generated by $\TyptrE_i$ (resp. $\TyptrF_i$) ($i\in\TypfkI$).
Define the $\bK$-linear homomorphism $\Typtrmtwo:\TypUm\otimes_\bK\TypUo\otimes_\bK\TypUp\to\TypUchipi$
by $\Typtrmtwo(Y\otimes Z\otimes X):=YZX$.
Then $\Typtrmtwo$ is a $\bK$-linear isomorphism.
\newline\newline
$(\TypU 6)$  For $\lambda\in\TypfkAoriginalpi$,
define $\TypUp_\lambda=\TypUpchipi_\lambda$
(resp. $\TypUm_\lambda=\TypUmchipi_\lambda$)
by $\TypUp_\lambda:=\TypUp\cap\TypU_\lambda$
(resp. $\TypUm_\lambda=\TypUm\cap\TypU_\lambda$).
Then for $\lambda\in\TypfkAoriginalpip\setminus\{0\}$,
we have 
$\{X\in\TypUp_\lambda|\forall i\in\TypfkI, [X,\TyptrF_i]=0\}=\{0\}$
and $\{Y\in\TypUm_{-\lambda}|\forall i\in\TypfkI, [\TyptrE_i,Y]=0\}=\{0\}$.
\newline\newline
Note that $\TypUp=\oplus_{\lambda\in\TypfkAoriginalpip}\TypUp_\lambda$,
$\TypUm=\oplus_{\lambda\in\TypfkAoriginalpip}\TypUm_{-\lambda}$,
$\TypUp_0=\TypUm_0=\bK\cdot1_\TypU$
and $\TypUp_{\Typpial_i}=\bK\cdot\TyptrE_i$,
$\TypUm_{-\Typpial_i}=\bK\cdot\TyptrF_i$
($i\in\TypfkI$).
\newline\par
We also regard $\TypU=\TypUchipi$ as a Hopf algebra $(\TypU,\TypDelta,\TypS,\Type)$
by
\begin{equation}\label{eqn:defHopf}
\begin{array}{l}
\TypDelta(\TyptrK_\lambda)=\TyptrK_\lambda\otimes\TyptrK_\lambda,
\TypDelta(\TyptrL_\lambda)=\TyptrL_\lambda\otimes\TyptrL_\lambda,
\TypDelta(\TyptrE_i)=\TyptrE_i\otimes  1+\TyptrK_{\Typpial_i}\otimes\TyptrE_i, \\
\TypDelta(\TyptrF_i)=\TyptrF_i\otimes
\TyptrL_{\Typpial_i}+1\otimes  \TyptrF_i,
\TypS(\TyptrK_\lambda)=\TyptrK_{-\lambda}, 
\TypS(\TyptrL_\lambda)=\TyptrL_{-\lambda}, \\
\TypS(\TyptrE_i)=-\TyptrK_{-\Typpial_i}\TyptrE_i,
\TypS(\TyptrF_i)=-\TyptrF_i\TyptrL_{-\Typpial_i}, \\
\Type(\TyptrK_\lambda)=\Type(\TyptrL_\lambda)=1,
\Type(\TyptrE_i)=\Type(\TyptrF_i)=0.
\end{array}
\end{equation}

Let $\TypUpflat:=\oplus_{\lambda\in\TypfkAoriginal}\TypUp\TyptrK_\lambda$,
and $\TypUmflat:=\oplus_{\lambda\in\TypfkAoriginal}\TypUm\TyptrL_\lambda$.
Then $\TypU=\TyprmSpan_\bK(\TypUmflat\TypUpflat)=\TyprmSpan_\bK(\TypUpflat\TypUmflat)$.

As in a
standard way (see \cite{Dr86}), we have the bilinear form
$\Typvartheta=\Typvarthetachipi:\TypUpflat\times\TypUmflat\to\bK$ 
such that $\Typvartheta_{|\TypUp\times\TypUm}$ is non-degenerate and
the following
equations are satisfied.
\begin{equation*}
\begin{array}{l}
\Typvartheta(\TyptrK_\lambda,\TyptrL_\mu)=\Typchi(\lambda,\mu),
\Typvartheta(\TyptrE_i,\TyptrF_j)=\delta_{ij},
\Typvartheta(\TyptrK_\lambda,\TyptrF_j)=\Typvartheta(\TyptrE_i,\TyptrL_\lambda)=0, \\
\Typvartheta(X^+Y^+,X^-)=
\sum_{k^-}\Typvartheta(X^+,(X^-)^{(2)}_{k^-})\Typvartheta(Y^+,(X^-)^{(1)}_{k^-}),\\
\Typvartheta(X^+,X^-Y^-)=
\sum_{k^+}\Typvartheta((X^+)^{(1)}_{k^+},X^-)\Typvartheta((X^+)^{(2)}_{k^+},Y^-),\\
\Typvartheta(\TypS(X^+),X^-)=\Typvartheta(X^+,\TypS^{-1}(X^-)),\\
\Typvartheta(X^+,1)=\Type(X^+),
\Typvartheta(1,X^-)=\Type(X^-), \\
X^-X^+  =\sum_{r^+,r^-}
\Typvartheta((X^+)^{\prime,(1)}_{r^+}, \TypS((X^-)^{\prime,(1)}_{r^-}))
\Typvartheta((X^+)^{\prime,(3)}_{r^+}, (X^-)^{\prime,(3)}_{r^-}) \\
\quad\quad\quad\quad\quad\quad\quad\quad
\cdot(X^+)^{\prime,(2)}_{r^+}(X^-)^{\prime,(2)}_{r^-}, \\
X^+X^- =\sum_{r^+,r^-}
\Typvartheta((X^+)^{\prime,(3)}_{r^+}, \TypS((X^-)^{\prime,(3)}_{r^-}))
\Typvartheta((X^+)^{\prime,(1)}_{r^+}, (X^-)^{\prime,(1)}_{r^-}) \\
\quad\quad\quad\quad\quad\quad\quad\quad
\cdot(X^-)^{\prime,(2)}_{r^-}(X^+)^{\prime,(2)}_{r^+}
\end{array}
\end{equation*} $(\lambda$, $\mu\in\TypfkAoriginal$, $i$, $j\in\TypfkI$,
and $X^+$, $Y^+\in\TypUpflat$, $X^-$, $Y^-\in\TypUmflat)$.
Here $(X^+)^{(x)}_{k^+}$ and $(X^-)^{(x)}_{k^-}$ with $x\in\TypfkJ_{1,2}$
(resp. $(X^+)^{\prime,(y)}_{r^+}$ and
$(X^-)^{\prime,(y)}_{r^-}$ with $y\in\TypfkJ_{1,3}$)
are any elements of $\TypUpflat$ and $\TypUmflat$ respectively
satisfying $\TypDelta(X^\pm)=\sum_{k^\pm}(X^\pm)^{(1)}_{k^\pm}
\otimes (X^\pm)^{(2)}_{k^\pm}$,
(resp. $((\Typrmid_\TypU\otimes\TypDelta)\circ\TypDelta)(X^\pm)=
\sum_{r^\pm}(X^\pm)^{\prime,(1)}_{r^\pm}\otimes
(X^\pm)^{\prime,(2)}_{r^\pm}\otimes (X^\pm)^{\prime,(3)}_{r^\pm}$).
We have
$\Typvartheta(X^+\TyptrK_\lambda,X^-\TyptrL_\mu)=\Typchi(\lambda,\mu)\Typvartheta(X^+,X^-)$
$(\lambda,\mu\in\TypfkAoriginal,\,X^+\in\TypUp,\,X^-\in\TypUm)$
and
$\Typvartheta(\TypUp_\lambda,\TypUm_{-\mu})=\{0\}$ $(\lambda.\mu\in\TypfkAoriginalpip$, $\lambda\ne\mu)$.

Define the $\bK$-linear map 
$\TypShchipi:\TypUchipi\to\TypUochipi$
by 
\begin{equation*}
\TypShchipi(Y\TyptrK_\lambda\TyptrL_\mu X)
=\Type(Y)\Type(X)\TyptrK_\lambda
\TyptrL_\mu
\end{equation*} 
$(X\in\TypUp,\,Y\in\TypUm,\,\lambda,\,\mu\in\TypfkAoriginal)$.

Let $\TypfkI^\prime$ be a non-empty subset of $\TypfkI$.
Let $\Typl_1:=|\TypfkI^\prime|$.
Define  the bijection $\Typkappa_{\TypfkI^\prime}:\TypfkJ_{1,\Typl_1}\to
\TypfkI^\prime$
by $\Typkappa_{\TypfkI^\prime}(t)<\Typkappa_{\TypfkI^\prime}(t+1)$
$(t\in\TypfkJ_{1,\Typl_1-1})$.
Let $\Typl_2\in\TypfkJ_{0,\Typlp-\Typl}$. If $\Typl_2\ne 0$, let 
$j_y\in\TypfkJ_{1,\Typlp-\Typl}$ $(y\in\TypfkJ_{1,\Typl_2})$ be such that 
$j_y<j_{y+1}$ $(y\in\TypfkJ_{1,\Typl_2-1})$.
Let $\TypfkAoriginal^\prime:=(\oplus_{t=1}^{\Typl_1}\bZ\Typpial_{\Typkappa_{\TypfkI^\prime}(t)})
\oplus(\oplus_{y=1}^{\Typl_2}\bZ\Typepsilon_{j_y})$.
Define the map $\Typpi_{\TypfkI^\prime}:\TypfkJ_{1,\Typl_1}\to\TypfkAoriginal^\prime$
by $\Typpi_{\TypfkI^\prime}:=\Typpi\circ\Typkappa_{\TypfkI^\prime}$.
Then we have the Hopf algebra monomorphism 
$f:\TypU(\Typchi_{|\TypfkAoriginal^\prime\times\TypfkAoriginal^\prime},
\Typpi_{\TypfkI^\prime})\to\TypUchipi$ defined by
$f(\TyptrK_\lambda\TyptrL_\mu):=\TyptrK_\lambda\TyptrL_\mu$
$(\lambda,\mu\in\TypfkAoriginal^\prime)$,
$f(\TyptrE_t):=\TyptrE_{\Typkappa_{\TypfkI^\prime}(t)}$, 
$f(\TyptrF_t):=\TyptrF_{\Typkappa_{\TypfkI^\prime}(t)}$
$(t\in\TypfkJ_{1,\Typl_1})$. So
\begin{equation}\label{eqn:ident}
\mbox{we identify $\TypU(\Typchi_{|\TypfkAoriginal^\prime\times\TypfkAoriginal^\prime},
\Typpi_{\TypfkI^\prime})$ with ${\mathrm{Im}}f$.}
\end{equation}

Let $\Typq_{ij}:=\Typchi(\Typpial_i,\Typpial_j)$
$(i,j\in\TypfkI)$.
We call $(\Typchi,\Typpi)$ {\it{reducible}}
if there exists non-empty proper subsets $\TypfkI^\prime$,
$\TypfkI^{\prime\prime}$
such that $\TypfkI^\prime\cap\TypfkI^{\prime\prime}=\emptyset$,
$\TypfkI=\TypfkI^\prime\cup\TypfkI^{\prime\prime}$ and
$\Typq_{ij}\Typq_{ji}=1$
for all $i\in\TypfkI^\prime$ and all $j\in\TypfkI^{\prime\prime}$. 
If $(\Typchi,\Typpi)$ is not reducible,
we call it {\it{irreducible}}. See also the map $f_3$ of Subsection~\ref{subsection:Gen}.

\subsection{Kharchenko's PBW theorem}

For $\lambda\in\TypfkAoriginal$, let $\Typq_\lambda:=\Typchi(\lambda,\lambda)$
and $\Typc_\lambda:=\Typkpch(\Typq_\lambda)$.

\begin{theorem}{\rm{(}}Kharchenko's PBW theorem~{\rm{\cite[Theorem~2]{Kha99}}}, {\rm{\cite[Theorem~2.2]{Kh15}}},
see also {\rm{\cite[Theorem~3.14]{HLec08}, \cite[Appendix]{BY19}}}.{\rm{)}} \label{theorem:KhPBW} 
Keep the notation as above. Then there exists a unique pair of $(\TypRchipip,\Typvphpip)$
of a subset $\TypRchipip$ of $\TypfkAoriginalpip\setminus\{0\}$ and a map $\Typvphpip:\TypRchipip\to\bN$ satisfying the following. 
Let $X:=\{(\al,t)\in\TypRchipip\times\bN|t\in\TypfkJ_{1,\Typvphpip(\al)}\}$.
Define the map $z:X\to\TypRchipip$ by $z(\al,t):=\al$.
Let $Y$ be the set of maps $y:X\to\bZgeqo$
such that $|\{x\in X|y(x)\geq 1\}|<\infty$ and $(y(x))_{\Typq_{z(x)}}!\ne 0$
for all $x\in X$. Then
\begin{equation*}
\forall\lambda\in\TypfkAoriginalpip,\,\dim\TypUpchipi_\lambda=|\{y\in Y|\sum_{x\in X}y(x)z(x)=\lambda\}|.
\end{equation*}
\end{theorem}

By \cite{Hec09}(see also \cite[Theorem~4.15~(1)]{AYY15}), using \eqref{eqn:elsRp} below, we see that
\begin{equation}\label{eqn:qbetak}
\begin{array}{l}
\mbox{if $\TyprmChar(\bK)=0$, $\Typl\geq 2$ and $|\TypRchipip|<\infty$
and $(\Typchi,\Typpi)$ is irreducible,} \\ 
\mbox{then $\Typqbeta\ne 1$ 
for all $\beta\in\TypRchipip$.}
\end{array}
\end{equation}

\begin{theorem}\label{theorem:PBWfinite}
{\rm{(\cite[Proposition~1]{Hec06}, \cite[Theorem~4.9]{HY10}}}{\rm{)}}
Assume $|\TypRchipip|<\infty$. \newline
{\rm{(1)}} We have $\Typvphpip(\TypRchipip)=\{1\}$.
\newline
{\rm{(2)}} Let $M:=|\TypRchipip|$.
Then there exist $\TyptrE_\beta\in\TypUpchipi_\beta$, $\TyptrF_\beta\in\TypUmchipi_{-\beta}$
$(\beta\in\TypRchipip)$
such that $\TyptrE_\beta\TyptrF_\beta-\TyptrF_\beta\TyptrE_\beta=-\TyptrK_\beta+\TyptrL_\beta$
and $\{\TyptrE_{g(1)}^{m_1}\cdots \TyptrE_{g(M)}^{m_M}| (m_s)_{\dummyTypq{g(s)}}!\ne 0\,(s\in\TypfkJ_{1,M})\}$
and $\{\TyptrF_{g(1)}^{n_1}\cdots \TyptrE_{g(M)}^{n_M}| (n_t)_{\dummyTypq{g(s)}}!\ne 0\,(t\in\TypfkJ_{1,M})\}$
are $\bK$-bases of $\TypUpchipi$ and $\TypUmchipi$ respectively
for every bijection $g:\TypfkJ_{1,M}\to\TypRchipip$.
\end{theorem}

Let  $\Typomega:\TypfkAoriginal\to\bKt$ be a $\bZ$-module homomorphism. 
Let
\begin{equation}\label{eqn:defprtrZomega}
\TypprtrZomchipi:=\{\,Z\in\TypUchipi_0\,|\,\forall \lambda\in\TypfkAoriginalpi,\,\forall X\in\TypUchipi_\lambda,\,
ZX
=\Typomega(\lambda)XZ\,\}.
\end{equation}
Define the $\bZ$-module homomorphism $\Typhrhochipi:\TypfkAoriginalpi\to\bKt$ by
\begin{equation*}
\Typhrhochipi(\Typpial_j):=\Typq_{\Typpial_j}\quad(j\in\TypfkI),
\end{equation*} where $\Typpial_j:=\Typpi(j)$, as above.
Let $\Typomegachhilmb
:=\Typomega(\beta)\cdot{\frac {\Typchi(\beta,\mu)} {\Typchi(\lambda,\beta)}}$
$(\beta\in\TypRchipip$, $\lambda$, $\mu\in\TypfkAoriginal)$.

For each $\beta\in\TypRchipip$, let $\TypfkBchipiomega(\beta)$ be the $\bK$-linear
subspace of $\TypUochipi$ formed by
the elements 
\begin{equation*}
\sum_{(\lambda,\mu)\in\TypfkAoriginal^2}a_{(\lambda,\mu)}
\TyptrK_\lambda\TyptrL_\mu
\end{equation*} with $a_{(\lambda,\mu)}\in\bK$
satisfying the following equations $(e1)_\beta$-$(e4)_\beta$.
\newline\newline
$(e1)_\beta$ For $(\lambda,\mu)\in\TypfkAoriginal^2$ and $t\in\bZ\setminus\{0\}$,
if $\Typqbeta\ne 1$,
$\Typcbeta=0$ 
and $\Typomegachhilmb=\Typqbeta^t$,
then the equation
$a_{(\lambda+t\beta,\mu-t\beta)}
=\Typhrhochipibeta^t\cdot a_{(\lambda,\mu)}$ holds.
\newline\newline
$(e2)_\beta$ For $(\lambda,\mu)\in\TypfkAoriginal^2$, if $\Typcbeta=0$ and 
$\Typomegachhilmb\ne\Typqbeta^t$ for all $t\in\bZ$, 
then the equation $a_{(\lambda,\mu)}=0$ holds.
\newline\newline
$(e3)_\beta$ For $(\lambda,\mu)\in\TypfkAoriginal^2$,
if $\Typqbeta\ne 1$, $\Typcbeta\geq 2$ and 
$\Typomegachhilmb=\Typqbeta^t$ for some $t\in\TypfkJ_{1,\Typcbeta-1}$, 
the equation
\begin{equation*}
\begin{array}{l}
\displaystyle{\sum_{x=-\infty}^\infty}a_{(\lambda+(\Typcbeta x+t)\beta,\mu-(\Typcbeta x+t)\beta)}
\Typhrhochipibeta^{-(\Typcbeta x+t)} \\
\quad =\displaystyle{\sum_{y=-\infty}^\infty}a_{(\lambda+\Typcbeta y\beta,\mu-\Typcbeta y\beta)}
\Typhrhochipibeta^{-\Typcbeta y}
\end{array}
\end{equation*} holds.
\newline\newline
$(e4)_\beta$ For $(\lambda,\mu)\in\TypfkAoriginal^2$, 
if  $\Typcbeta\geq 2$ and 
$\Typomegachhilmb\ne\Typqbeta^m$ for all $m\in\bZ$, then the $\Typcbeta-1$ equations
\begin{equation*}
\begin{array}{l}
\displaystyle{\sum_{x=-\infty}^\infty}a_{(\lambda+(\Typcbeta x+t)\beta,\mu-(\Typcbeta x+t)\beta)}
\Typhrhochipibeta^{-(\Typcbeta x+t)} \\
\quad =\displaystyle{\sum_{y=-\infty}^\infty}a_{(\lambda+\Typcbeta y\beta,\mu-\Typcbeta y\beta)}
\Typhrhochipibeta^{-\Typcbeta y} \\
(t\in\TypfkJ_{1,\Typcbeta-1})
\end{array}
\end{equation*} hold.
Let
\begin{equation*}
\TypfkBchipiomega:=\bigcap_{\beta\in\TypRchipip}\TypfkBchipiomega(\beta).
\end{equation*}

\begin{theorem}\label{theorem:previousmain}
{\rm{(}}{\rm{\cite[Theorem~10.4]{BY18}}}{\rm{)}}
Assume $|\TypRchipip|<\infty$.
If ${\mathrm{Char}}(\bK)\ne 0$, assume that $\Typqbeta\ne 1$ for all $\beta\in\TypRchipip$.
Then we have the $\bK$-linear isomorphism $\TypHCbhmpiomega:\TypprtrZomchipi\to\TypfkBchipiomega$ defined by $\TypHCbhmpiomega(X):=\TypShchipi(X)$.
\end{theorem}
{\it{Proof.}} 
Injectivity of $\TypHCbhmpiomega$ can be proved in the same way as that for \cite[Lemma~9.2]{BY18}.

First we assume that $\TypfkAoriginalEX=\{0\}$ and  $\Typqbeta\ne 1$ for all $\beta\in\TypRchipip$.
The statement of \cite[Theorem~10.4]{BY18} has claimed that the above theorem holds if $\bK$ is an algebraically closed field.
However, by the equation ${\mathcal{G}}={\mathcal{Z}}{\mathcal{S}}$ in \cite[(10.5)]{BY18},
we can easily see that it really holds for any field $\bK$.

Second we assume $\TypfkAoriginalEX=\{0\}$ and  ${\mathrm{Char}}(\bK)=0$.
We use an induction on $|\TypfkI|$.
Assume that there exists $\beta\in\TypRchipip$ with $\Typqbeta=1$.
By \eqref{eqn:qbetak}, there exists $i\in\TypfkI$ such that
$\beta=\Typpial_i$ 
and $\Typq_{ij}\Typq_{ji}=1$ for all $j\in\TypfkI\setminus\{i\}$.
Let ${\acute{\TypfkI}}:=\TypfkI\setminus\{i\}$.
Then $\TyptrE_i\TyptrE_j=\Typq_{ij}\TyptrE_j\TyptrE_i={\frac 1 {\Typq_{ji}}}\TyptrE_j\TyptrE_i$
and $\TyptrF_i\TyptrF_j=\Typq_{ji}\TyptrF_j\TyptrF_i={\frac 1 {\Typq_{ij}}}\TyptrF_j\TyptrF_i$ for all $j\in{\acute{\TypfkI}}$.
Let ${\grave{\TypU}}_0:=\oplus_{k=0}^\infty\TyptrF_i^k\TypUochipi\TyptrE_i^k$.
Let 
${\acute{\mathfrak{A}}}_\Typpi:=\oplus_{j\in{\acute{\TypfkI}}}\bZ\Typpial_j$.
${\acute{\mathfrak{A}}}_\Typpi^+:=\oplus_{j\in{\acute{\TypfkI}}}\bZgeqo\Typpial_j$,
${\acute{\Typchi}}:=\Typchi_{|{\acute{\mathfrak{A}}}_\Typpi\times{\acute{\mathfrak{A}}}_\Typpi}$,
and ${\acute{\Typpi}}:=\Typpi_{{\acute{\TypfkI}}}$.
Let ${\acute{\TypU}}_0^\lambda:=\TyprmSpan_\bK(\TypUmchipi_{-\lambda}\TypUpchipi_\lambda)$
$(\lambda\in{\acute{\mathfrak{A}}}_\Typpi^+)$. Let
${\acute{\TypU}}_0:=\oplus_{\lambda\in{\acute{\mathfrak{A}}}_\Typpi^+}{\acute{\TypU}}_0^\lambda$.
Then 
\begin{equation}\label{eqn:facTypU}
\begin{array}{l}
\mbox{we have the $\bK$-linear isomorphism} \\
\mbox{$f:{\acute{\TypU}}_0\otimes{\grave{\TypU}}_0\to\TypUchipi_0$
defined by $f(X\otimes Y):=XY$.}
\end{array}
\end{equation}
Let $a\in\bKt$. Let 
${\dot{\mathfrak{Z}}}_a:=\{\,Z\in\TypUchipi_0\,|
\,ZE_i=aE_iZ,\,ZF_i=a^{-1}F_iZ\}$,
and ${\grave{\mathfrak{Z}}}_a:=\{C\in{\grave{\TypU}}_0|\,CE_i=aE_iC,\,CF_i=a^{-1}F_iC\}$.
Let ${\acute{\mathcal{H}}}_a:=\oplus_{\lambda,\mu\in{\acute{\mathfrak{A}}}_\Typpi, 
a{\frac{\Typchi(\Typpial_i,\mu)} {\Typchi(\lambda,\Typpial_i)}}=1}
\bK\TyptrK_\lambda\TyptrL_\mu$,
and ${\grave{\mathcal{H}}}_a:=\oplus_{\lambda,\mu\in\TypfkAoriginalpi,
a{\frac{\Typchi(\Typpial_i,\mu)} {\Typchi(\lambda,\Typpial_i)}}=1}
\bK\TyptrK_\lambda\TyptrL_\mu$.
Then ${\grave{\mathcal{H}}}_a=\oplus_{x,y\in\bZ}{\acute{\mathcal{H}}}_a
\TyptrK_{\Typpial_i}^x\TyptrL_{\Typpial_i}^y$.
By \cite[Theorem~2.4]{BY15}, we have 
${\grave{\mathfrak{Z}}}_a={\grave{\mathcal{H}}}_a$.
By \eqref{eqn:facTypU}, we have ${\dot{\mathfrak{Z}}}_a=
\oplus_{\Typchi(\lambda,\Typpial_i)b=a}
\TyprmSpan_\bK({\acute{\TypU}}_0^\lambda{\grave{\mathcal{H}}}_b)$,
which implies $\TypShchipi({\dot{\mathfrak{Z}}}_a)={\grave{\mathcal{H}}}_a$.
Since $\TypprtrZomchipi\subset{\dot{\mathfrak{Z}}}_{\Typomega(\Typpial_i)}$, 
we have
$\TypShchipi(\TypprtrZomchipi)
\subset{\grave{\mathcal{H}}}_{\Typomega(\Typpial_i)}$.
The condition $(e2)_{\Typpial_i}$ means that $a_{(\lambda,\mu)}=0$
if $\Typomega(\Typpial_i)\cdot{\frac {\Typchi(\Typpial_i,\mu)} {\Typchi(\lambda,\Typpial_i)}}\ne 1$.
Hence ${\grave{\mathcal{H}}}_{\Typomega(\Typpial_i)}=\{Z\in\TypUochipi|\mbox{$Z$
satisfies $(e2)_{\Typpial_i}$}\}$.
For $x$, $y\in\bZ$, define the group homomorphism ${\acute{\Typomega}}_{x,y}:{\acute{\mathfrak{A}}}_\Typpi\to\bKt$
by ${\acute{\Typomega}}_{x,y}(\lambda):=\Typchi(\Typpial_i,\lambda)^{-(x+y)}\Typomega(\lambda)$
$(\lambda\in{\acute{\mathfrak{A}}}_\Typpi)$. Note
that ${\acute{\Typomega}}_{x,y}(\lambda)={\frac {\Typchi(\lambda,y\Typpial_i)} {\Typchi(x\Typpial_i,\lambda)}}\Typomega(\lambda)$ $(\lambda\in{\acute{\mathfrak{A}}}_\Typpi)$.
Then 
\begin{equation*}
\TypfkBchipiomega=\oplus_{x,y\in\bZ}({\mathfrak{B}}^\acutechipi_{{\acute{\Typomega}}_{x,y}}
\cap{\acute{\mathcal{H}}}_{\Typomega(\Typpial_i)})
\TyptrK_{\Typpial_i}^x\TyptrL_{\Typpial_i}^y.
\end{equation*}
Since $\TypprtrZomchipi$
$\subset\oplus_{x,y\in\bZ}
{\mathfrak{Z}}_{{\acute{\Typomega}}_{x,y}}(\acutechipi)\TyptrK_{\Typpial_i}^x\TyptrL_{\Typpial_i}^y$,
we have $\TypShchipi(\TypprtrZomchipi)\subset\TypfkBchipiomega$.
By \cite[Theorem~10.4]{BY18}  and a careful understanding 
(see Remark~\ref{remark:GeqZS} below) 
the equation ${\mathcal{G}}={\mathcal{Z}}{\mathcal{S}}$ in \cite[(10.5)]{BY18}, 
we see that
for every $Z\in{\mathfrak{B}}^\acutechipi_{{\acute{\Typomega}}_{x,y}}
\cap{\acute{\mathcal{H}}}_{\Typomega(\Typpial_i)}$, there exists 
$C\in{\mathfrak{Z}}_{{\acute{\Typomega}}_{x,y}}(\acutechipi)$
such that $\TypShchipi(C)=Z$ and $C\TyptrE_i=\Typomega(\Typpial_i)E_iC$
and $C\TyptrF_i={\frac 1 {\Typomega(\Typpial_i)}}\TyptrF_iC$.
Then we can see that the statement holds. 

Finally we assume $\TypfkAoriginalEX\ne \{0\}$.
For $\lambda^\prime$, $\mu^\prime\in\TypfkAoriginalEX$,
define the $\bZ$-module homomorphism 
$\Typomega_{\lambda^\prime,\mu^\prime}:\TypfkAoriginalpi\to\bKt$
by $\Typomega_{\lambda^\prime,\mu^\prime}(\lambda):=\Typomega(\lambda)
{\frac{\Typchi(\lambda,\mu^\prime)} {\Typchi(\lambda^\prime,\lambda)}}$. 
Let $\Typchi^{\prime\prime}:=\Typchi_{|\TypfkAoriginalpi\times\TypfkAoriginalpi}$.
Since $\TypUchipi=\oplus_{\lambda^\prime, \mu^\prime\in\TypfkAoriginalEX}
\TyptrK_{\lambda^\prime}\TyptrL_{\mu^\prime}\TypU(\Typchi^{\prime\prime},\Typpi)$,
we have
\begin{equation*}
\TypprtrZomchipi=\oplus_{\lambda^\prime, \mu^\prime\in\TypfkAoriginalEX}\TyptrK_{\lambda^\prime}\TyptrL_{\mu^\prime}\TypprtrZ_{\Typomega_{\lambda^\prime,\mu^\prime}}
(\Typchi^{\prime\prime},\Typpi).
\end{equation*}
Then we can see that the statement is true.
\hfill $\Box$

\begin{remark}\label{remark:GeqZS}
Seeing carefully \cite[Subsection~10.1]{BY18}, we see that for ${\mathcal{G}}={\mathcal{Z}}{\mathcal{S}}$ in \cite[(10.5)]{BY18}
for $\lambda\in{\acute{\mathfrak{A}}}_\Typpi$, 
the components of the matrix ${\mathcal{G}}$
(resp. ${\mathcal{Z}}$, resp. ${\mathcal{S}}$) belong to 
${\acute{\mathcal{H}}}_a$ for $a=\Typomega(\Typpial_i)$
(resp. ${\frac {\Typomega(\Typpial_i)} {\Typchi(\lambda,\Typpial_i)}}$, 
resp. $\Typchi(\lambda,\Typpial_i)$).
\end{remark}

\subsection{Basic facts}\label{subsection:Bf}

Let $\bKtinf:=\{x\in\bKt|x\ne 1, \Typkpch(x)\ne 0\}
=\{x\in\bKt|x^m\ne 1(m\in\bN)\}$.

In Subsection~\ref{subsection:Bf}, assume $|\TypRchipip|<\infty$.
Let $i\in\TypfkI$. For $j\in\TypfkI\setminus\{i\}$, let
$\TypNchipi_{ij}=\max\{m\in\bZgeqo|\Typpial_j+m\Typpial_i\in\TypRchipip\}$.
Define the map $\Typchitauipi:\TypfkI\to\TypfkAoriginalpi$
by $\Typchitauipi(i):=-\al_i$ and $\Typchitauipi(j):=\Typpial_j+\TypNchipi_{ij}\Typpial_i$
$(i\in\TypfkI\setminus\{j\})$. Clearly $\Typchitauipi(\TypfkI)$ is a $\bZ$-base of $\TypfkAoriginalpi$.
It is well-known that
\begin{equation}\label{eqn:tautaupi}
\begin{array}{l}
\mbox{$\Typchitau_i(\Typchitauipi)=\Typpi$, 
$\TypRchipiptaui\setminus\{-\Typpial_i\}=\TypRchipip\setminus\{\Typpial_i\}$,} \\
\mbox{$\TypNchipi_{ij}=\TypNchipitaui_{ij}=\max\{m\in\bZgeqo|(m)_{\Typq_{ii}}!(m;\Typq_{ii},\Typq_{ij}\Typq_{ji})!\ne 0\}$
$(j\in\TypfkI\setminus\{i\})$},
\end{array}
\end{equation}
see \cite[Proposition~1]{Hec06} (and \cite[Lemma~4.9~(2)]{AYY15}).
\begin{lemma}\label{lemma:fctGA}
{\rm{(}}See {\rm{\cite[Lemmas~4.21 and 4.22]{AYY15}}} for example.{\rm{)}}
Let $i\in\TypfkI$. Assume that $\Typq_{ii}\in\bKtinf$.
Let $\Typpial^\prime_j:=\Typchitauipi(j)$ $(j\in\TypfkI)$.
Let $\Typq^\prime_{jk}:=\Typchi(\Typpial^\prime_j,\Typpial^\prime_k)$
$(j,k\in\TypfkI)$. 
Define the $\bZ$-module isomorphism $\Typacsi:\TypfkAoriginalpi\to\TypfkAoriginalpi$
by $\Typacsi(\Typpial_j):=\Typpial^\prime_j$ $(j\in\TypfkI)$.
Then we have the following.
\newline 
{\rm{(1)}} For every $\lambda\in\TypfkAoriginalpi$, there exists $n\in\bZ$
such that $\Typchi(\lambda,\al_i)\Typchi(\al_i,\lambda)=\Typq_{ii}^n$.
\newline 
{\rm{(2)}} We have $q_{jj}=q^\prime_{jj}$ and $q_{jk}q_{kj}=q^\prime_{jk}q^\prime_{kj}$
for all $j$, $k\in\TypfkI$. \newline
{\rm{(3)}} We have $\Typacsi^2=\Typrmid_{\TypfkAoriginalpi}$, 
and 
\begin{equation*}
\Typchi(\Typacsi(\lambda),\Typacsi(\lambda))=\Typchi(\lambda,\lambda),\,\,
\Typchi(\Typacsi(\lambda),\Typacsi(\mu))\Typchi(\Typacsi(\mu),\Typacsi(\lambda))=
\Typchi(\lambda,\mu)\Typchi(\mu,\lambda)
\end{equation*}
for all $\lambda$, $\mu\in\TypfkAoriginalpi$.
\newline
{\rm{(4)}} We have 
$\Typacsi(\TypRchipip)=\TypRchipiptaui$.
Moreover $\Typacsi(\lambda)=\lambda-n\Typpial_i$
for $\lambda\in\TypfkAoriginalpi$ with
$\Typchi(\Typpial_i,\lambda)\Typchi(\lambda,\Typpial_i)=\Typq_{ii}^n$
and $n\in\bZ$.
\end{lemma}

Let $i\in\TypfkI$.
For $j\in\TypfkI\setminus\{i\}$, let 
$a_j$, $b_j\in\bKt$ be such that
\begin{equation*}
a_jb_j(\TypNchipi_{ij})_{\Typq_{ii}}!(\TypNchipi_{ij};\Typq_{ii},\Typq_{ij}\Typq_{ji})!=1,
\end{equation*}
and define $\TyptrEchipi_{m,i,j}\in\TypUpchipi_{\Typpial_j+m\Typpial_i}$
(resp. $\TyptrFchipi_{m,i,j}\in\TypUpchipi_{-\Typpial_j-m\Typpial_i}$)
$(m\in\bZgeqo)$ by
$\TyptrEchipi_{0,i,j}:=\TyptrE_j$
(resp. $\TyptrFchipi_{0,i,j}:=\TyptrF_j$)
and $\TyptrEchipi_{m+1,i,j}:=\TyptrE_i\TyptrEchipi_{m,i,j}-\Typq_{ii}^m\Typq_{ij}\TyptrEchipi_{m,i,j}\TyptrE_i$,
(resp. $\TyptrFchipi_{m+1,i,j}:=\TyptrF_i\TyptrFchipi_{m,i,j}-\Typq_{ii}^m\Typq_{ij}\TyptrFchipi_{m,i,j}\TyptrF_i$).
Then we have the $\bK$-algebra isomorphism
$\TypTchitauipii:\TypUchipitaui\to\TypUchipi$
with $\TypTchitauipii(\TyptrK_\lambda\TyptrL_\mu)=\TyptrK_\lambda\TyptrL_\mu$
$(\lambda,\mu\in\TypfkAoriginal)$,
$\TypTchitauipii(\TyptrE_i)=\TyptrF_i\TyptrL_{-\Typpial_i}$,
$\TypTchitauipii(\TyptrF_i)=\TyptrK_{-\Typpial_i}\TyptrE_i$,
$\TypTchitauipii(\TyptrE_j)=a_j\TyptrEchipi_{\TypNchipi_{ij},i,j}$ $(j\ne i)$
and $\TypTchitauipii(\TyptrF_j)=b_j\TyptrFchipi_{\TypNchipi_{ij},i,j}$ $(j\ne i)$,
see \cite[Theorem~6.7]{Hec10}.

Let $\TypbBI$ the set of all maps from $\bN$ to $\TypfkI$.
Let $f\in\TypbBI$ and $t\in\bZgeqo$. Define 
$\TypchitauftTyppi$ in the way that 
$\TypchitaufoTyppi:=\Typpi$ 
and $\TypchitauftTyppi:=\Typchitau_{f(t)}(\TypchitauftmTyppi)$.

Let $M:=|\TypRchipip|$. We know that
\begin{equation}\label{eqn:elsRp}
\mbox{there exists $f\in\TypbBI$ such that $\TypRchipip
=\{\TypchitauftmTyppi(f(t))|t\in\TypfkJ_{1,M}\}$,}
\end{equation} see \cite[Lemma~5.5]{BY18}.

Let $\TypRchipi:=\TypRchipip\cup(-\TypRchipip)$.
Let $\TypRchipire:=\{\beta\in\TypRchipi|\Typqbeta\in\bKtinf\}$,
$\TypRchipinu:=\TypRchipi\setminus\TypRchipire$,
and
$\TypRchipipre:=\TypRchipip\cap\TypRchipire$,
$\TypRchipipnu:=\TypRchipip\cap\TypRchipinu$.

\section{Main Result}
\subsection{Irreducible bicharacters
with $\TyprmChar(\bK)=0$ and $|\TypRchipip|<\infty$}
\label{subsection:datum}
Let $(\TypRchipire)^\flat:=\TypRchipire\cup\{\beta\in\TypRchipinu|\Typqbeta=1\}$.
\begin{equation}\label{eqn:assumption-a}
\begin{array}{l}
\mbox{From now on until Remark~\ref{remark:rankonecase}, 
we assume that} \\
\mbox{$\TyprmChar(\bK)=0$, $|\TypRchipip|<\infty$,
$\dim\TypUpchipi=\infty$
$(\Leftrightarrow(\TypRchipire)^\flat\ne\emptyset)$,
and} \\ 
\mbox{$(\Typchi,\Typpi)$ is irreducible.} \\
\end{array}
\end{equation}
By the classification theorem of \cite{Hec09}, 
we have an injection $\Typpibar:\TypfkI\to\TypfkAoriginal$ being one of $(\Typpibar 0)$-$(\Typpibar 5)$
below
such that $\TypfkAoriginalpi=\oplus_{i\in I}\bZ\Typpibar(i)$ and  
$\Typpi=\Typtau^\Typchi_{f,t}(\Typpibar)\circ d$
for some bijection $d:\TypfkI\to\TypfkI$, some $f\in\TypbBI$
and some $t\in\bZgeqo$
after modifying $\Typlp$ and the values 
$\Typchi(\Typpi(i)), \Typepsilon_r)$,
$\Typchi(\Typepsilon_r,\Typpi(i))$,
$\Typchi(\Typepsilon_r,\Typepsilon_{r^\prime})$.
\begin{equation}\label{eqn:assumption-b}
\begin{array}{l}
\mbox{From now on until Remark~\ref{remark:rankonecase}, 
using $(\Typpibar 0)$-$(\Typpibar 5)$ below,} \\
\mbox{we fix $\Typlp$, $\Typchi(\Typpibar(i),\Typepsilon_r)$,
$\Typchi(\Typepsilon_r,\Typpibar(i))$,
$\Typchi(\Typepsilon_r,\Typepsilon_{r^\prime})$.
See also \eqref{eqn:relEX} below.}
\end{array}
\end{equation}
By \eqref{eqn:tautaupi}, we have
\begin{equation}\label{eqn:pibarpi}
\TypRchibarpi=\TypRchipi,\,\,\TypRchibarpire=\TypRchipire,\,\,\TypRchibarpinu=\TypRchipinu,\,\,
(\TypRchibarpi)^\flat=(\TypRchipi)^\flat.
\end{equation}

Let $\Typpialbar_i:=\Typpibar(i)$.
Let $\Typqbar_{ij}:=\Typchi(\Typpialbar_i,\Typpialbar_j)$.
Let $\TypfkAoriginalbR:=\TypfkAoriginal\otimes_\bZ\bR$.
Then $\TypfkAoriginalbR$ is an $\bR$-linear space.
Note that the $\bZ$-module homomorphism 
$\TypUpsilon:\TypfkAoriginal\to\TypfkAoriginalbR$ defined by $\TypUpsilon(\lambda):=\lambda\otimes 1$ 
$(\lambda\in\TypfkAoriginal)$ is injective.
For $\lambda\in\TypfkAoriginal$, we denote $\TypUpsilon(\lambda)$ by $\lambda$ for simplicity.
We define $\Typvarpi_i\in
\TypfkAoriginal$ $(i\in\TypfkI)$ and a non-degenerate symmetric bilinear map 
$(\,|\,):\TypfkAoriginalbR\times\TypfkAoriginalbR\to\bR$
by using $(\Typpibar 0)$-$(\Typpibar 5)$ below.
Let $\Typbark:=\det[(\Typpialbar_i|\Typpialbar_j)]_{i,j\in I}$.
\newline\newline
$(\Typpibar 0)$ ($\Typl=\Typqbar_{11}=1$-case)
Assume $\Typl = 1$ and $\Typqbar_{11}=1$. 
We have $\TypRchibarpipre=\emptyset$
and $\TypRchibarpipnu=\{\Typpialbar_1\}$.
Let $x\in\bKtinf$.
Let $\Typlp:=2$, $\Typvarpi_1:=\Typepsilon_1$,
$\TypfkAoriginalEX=\bZ\Typvarpi_1$, 
$\Typchi(\Typvarpi_1,\Typvarpi_1)=1$,
$\Typchi(\Typpialbar_1,\Typvarpi_1):=1$, $\Typchi(\Typvarpi_1,\Typpialbar_1)=x$,
and $(\Typvarpi_1|\Typvarpi_1)=(\Typpialbar_1|\Typpialbar_1)=0$
and $(\Typvarpi_1|\Typpialbar_1)=1$.
\newline\newline
$(\Typpibar 1)$ (Cartan-type)
Let $\Typlp:=\Typl$.
Define $(\,|\,)$
in the way that  $(\Typpialbar_j|\Typpialbar_j)\in\bN$ $(j\in\TypfkI)$ and
$A:=[{\frac {2(\Typpialbar_i|\Typpialbar_j)} {(\Typpialbar_j|\Typpialbar_j)}}]_{i,j\in\TypfkI}$ is one on the matrices in 
\cite[Table 1 of Subsection~1.4]{Hum}.
Namely $A$ is the Cartan matrix of an irreducible root system in the sense of \cite[Section~11]{Hum}.
Let $\Typchi$ be such that
$\Typqbar_{ii}=\Typq^{(\Typpialbar_i|\Typpialbar_i)}$ and
$\Typqbar_{ij}\Typqbar_{ji}=\Typq^{2(\Typpialbar_i|\Typpialbar_j)}$
$(i,j\in\TypfkI)$ for some $\Typq\in\bKtinf$. 
Note that $(\,|\,)$ is positive definite.
Let $\Typvarpi_i\in\TypfkAoriginal$ 
$(i\in\TypfkI)$ be such that 
$(\Typvarpi_i|\Typpialbar_j)=\delta_{ij}\cdot\Typbark$
$(i,j\in\TypfkI)$.
We have $\TypRchibarpipnu=\emptyset$.
The set $\TypRchibarpip=\TypRchibarpipre$ 
can be identified with the positive root system of 
the finite-dimensional simple Lie algebra over $\bC$
whose Cartan matrix is $A$,
see \cite[Lemma~5.3]{AYY15} for example.
\newline\newline
$(\Typpibar 2)$ (Super-type)
Let $X=[x_{ij}]_{ij\in\TypfkI}$ be one of the $\Typl\times\Typl$-symmetric matrices 
(i.e., $x_{ij}=x_{ji}$) over $\bZ$ 
below.
\newline\newline
{\rm{(i)}} (${\mathfrak{sl}}(m|n)$-type) Assume $\Typl\geq 2$. Let $m$, $n\in\bN$
be such that 
$m+n+1=\Typl$ and $m\geq n$.
Let $x_{ii}:=2$, $x_{ii+1}:=-1$ $(i\in\TypfkI_{1,m})$,
$x_{m+1,m+1}:=0$,
$x_{jj}:=-2$, $x_{j-1,j}:=1$ $(j\in\TypfkI_{m+2,m+n+1})$,
and $x_{i^\prime,j^\prime}:=0$
$(|i^\prime-j^\prime|\geq 2)$.
\newline
{\rm{(ii)}} ($B(m,n)$-type) Assume $\Typl\geq 2$. 
Let $m$, $n\in\bN$ be such that $m+n=\Typl$.
Let $x_{ii}:=-2$, $x_{ii+1}:=1$ $(i\in\TypfkI_{1,n-1})$,
$x_{nn}:=0$,
$x_{jj}:=2$, $x_{j-1,j}:=-1$ $(i\in\TypfkI_{n+1,m+n-1})$,
$x_{m+n,m+n}:=1$,  $x_{m+n-1,m+n}:=-1$
and $x_{i^\prime,j^\prime}:=0$
$(|i^\prime-j^\prime|\geq 2)$.
\newline
{\rm{(iii)}} ($C(n)$-type) Assume $\Typl\geq 3$. Let $n:=\Typl$.
Let $x_{11}:=0$,
$x_{ii}:=2$, $x_{i-1,i}:=-1$ $(j\in\TypfkI_{2,n-1})$,
$x_{nn}:=4$, $x_{n-1,n}:=-2$
and $x_{i^\prime,j^\prime}:=0$
$(|i^\prime-j^\prime|\geq 2)$.
\newline
{\rm{(iv)}} ($D(m,n)$-type) Assume $\Typl\geq 3$. 
Let $m\in\TypfkJ_{2,\infty}$, $n\in\bN$ be such that 
$m+n=\Typl$.
Let $x_{ii}:=-2$, $x_{ii+1}:=1$ $(i\in\TypfkI_{1,n-1})$,
$x_{nn}:=0$,
$x_{jj}:=2$, $x_{j-1,j}:=-1$ $(j\in\TypfkI_{n+1,m+n-1})$,
$x_{m+n,m+n}:=2$, $x_{m+n-1,m+n}:=0$, $x_{m+n-2,m+n}:=-1$, 
and $x_{i^\prime,j^\prime}:=0$
$(i^\prime+2\leq j^\prime,\, (i^\prime,j^\prime)\ne(m+n-2,m+n))$.
\newline
{\rm{(v)}} ($F(4)$-type) Assume $\Typl=4$.
Let $x_{11}:=0$,  $x_{22}:=2$,  $x_{33}:=x_{44}:=4$,  
$x_{12}:=-1$, $x_{23}:=x_{34}:=-2$,
and $x_{i^\prime,j^\prime}:=0$
$(|i^\prime-j^\prime|\geq 2)$.
\newline
{\rm{(vi)}} ($G(3)$-type) Assume $\Typl=3$.
Let $x_{11}:=0$,  $x_{22}:=2$,  $x_{33}:=6$,  
$x_{12}:=-1$, $x_{23}:=-3$,
and $x_{i^\prime,j^\prime}:=0$
$(|i^\prime-j^\prime|\geq 2)$.
\newline
{\rm{(vii)}}  ($D(2,1;a)$-type) Assume $\Typl=3$.
Let $a\in\bZ\setminus\{0,-1\}$.
Let $x_{11}:=-2$,  $x_{22}:=0$,  $x_{33}:=-2a$, 
$x_{12}:=1$, $x_{23}:=a$ and $x_{13}:=0$. 
\newline\newline
Note that $\det X= 0$ if and only if
$X$ is of ${\mathfrak{sl}}(m|m)$-type with $\Typl=2m+1$.
Define $(\,|\,)_{|\TypfkAoriginalpi\times\TypfkAoriginalpi}$ by
$(\Typpialbar_i|\Typpialbar_j):=x_{ij}$.
Let $\Typq\in\bKtinf$. 
Let $\Typchi_{|\TypfkAoriginalpi\times\TypfkAoriginalpi}$ be such that 
$\Typq_{ii}=\Typq^{x_{ii}}$ $(x_{ii}\ne 0)$,
$\Typq_{jj}=-1$ $(x_{jj}=0)$,
$\Typq_{i^\prime j^\prime}\Typq_{j^\prime i^\prime}=\Typq^{2x_{i^\prime j^\prime}}$. 
If $\det X\ne 0$,
let $\Typlp:=\Typl$,
and let $\Typvarpi_i\in\TypfkAoriginal$ 
$(i\in\TypfkI)$ such that 
$(\Typvarpi_i|\Typpialbar_j)=\delta_{ij}\cdot\Typbark$
$(i,j\in\TypfkI)$.
Assume $\det X=0$.
Let  $\Typlp=\Typl+1(=2m+2)$.
Let $(\Typepsilon_1|\Typepsilon_1):=1$ and $(\Typepsilon_1|\Typpialbar_i):=\delta_{i1}$ $(i\in\TypfkI)$.
Let $\Typchi(\Typepsilon_1,\Typepsilon_1):=\Typq$ and 
$\Typchi(\Typepsilon_1,\Typpialbar_i):=1$, $\Typchi(\Typpialbar_i,\Typepsilon_1):=\Typq^{2\delta_{i1}}$ $(i\in\TypfkI)$.
For $i\in\TypfkJ_{2,\Typl}$, let
$\Typepsilon_i:=\Typepsilon_1-\sum_{t=1}^{i-1}\Typpial_t$.
Let $\Typvarpi_i:=\sum_{t=1}^i\Typepsilon_t$  $(i\in\TypfkJ_{1,m+1})$
and $\Typvarpi_j:=\Typvarpi_{m+1}-\sum_{t=m+2}^i\Typepsilon_t$  $(i\in\TypfkJ_{m+2,\Typl})$.
We see that $(\Typvarpi_i|\Typpialbar_j)=\delta_{ij}$
and $\Typchi(\Typvarpi_i,\Typpialbar_j)\Typchi(\Typpialbar_j,\Typvarpi_i)=\Typq^{2\delta_{ij}}$.
Note that $\TypRchibarpi=\{\Typepsilon_i-\Typepsilon_j|i,j\in\TypfkI,i<j\}$. 
See also Remark~\ref{remark:relofInt} blow.
\newline\newline
$(\Typpibar 3)$ (Extra-$D(2,1;a)$-type) Let $x$, $y\in\bKtinf$ with $xy\ne 1$.
We assume that $\Typl=3$, $\Typq_{11}=x$,
$\Typq_{22}=-1$, $\Typq_{33}=y$,
$\Typq_{12}\Typq_{21}=x^{-1}$,
$\Typq_{23}\Typq_{32}=y^{-1}$,
$\Typq_{13}\Typq_{31}=1$.
Define $(\,|\,)_{|\TypfkAoriginalpi\times\TypfkAoriginalpi}$
by
$(\Typpialbar_1|\Typpialbar_1)=2$,
$(\Typpialbar_2|\Typpialbar_2)=0$,
$(\Typpialbar_3|\Typpialbar_3)=-4$,
$(\Typpialbar_1|\Typpialbar_2)=-1$,
$(\Typpialbar_1|\Typpialbar_3)=0$,
$(\Typpialbar_2|\Typpialbar_3)=2$. \newline
{\rm{(i)}} Assume  $xy\in\bKtinf$. Let 
$\Typlp=\Typl$, $\Typvarpi_1:=\Typpialbar_1$,
$\Typvarpi_2:=\Typpialbar_1+2\Typpialbar_2+\Typpialbar_3$,
$\Typvarpi_3:=\Typpialbar_3$.
Note $\Typchi(\Typvarpi_2,\Typpialbar_i)\Typchi(\Typpialbar_i,\Typvarpi_2)=(xy)^{-\delta_{i2}}$.
Note $\TypRchipipre=\{\Typvarpi_i|i\in\TypfkI\}$. 
We have $\TypRchipipnu=\{\Typpialbar_2,\Typpialbar_1+\Typpialbar_2,
\Typpialbar_2+\Typpialbar_3,\Typpialbar_1+\Typpialbar_2+\Typpialbar_3\}$. \newline
{\rm{(ii)}} Assume $xy\notin\bKtinf$. Let  $\Typlp:=\Typl+1$, $\Typvarpi_1:=\Typpialbar_1$,
$\Typvarpi_2:=\Typepsilon_1$,
$\Typvarpi_3:=\Typpialbar_3$,
$\Typchi(\Typepsilon_1,\Typepsilon_1):=1$,
$\Typchi(\Typpialbar_i,\Typepsilon_1):=1$, $\Typchi(\Typepsilon,\Typpialbar_i)=x^{\delta_{i2}}$,
and $(\Typepsilon_1|\Typepsilon_1)=(\Typepsilon_1|\Typpialbar_1)=(\Typepsilon_1|\Typpialbar_3)=0$
and $(\Typepsilon_1|\Typpialbar_2)=1$.
We have $\TypRchipipre=\{\Typpialbar_1,\Typpialbar_3\}$
and $\TypRchipipnu=\{\Typpialbar_2,\Typpialbar_1+\Typpialbar_2,
\Typpialbar_2+\Typpialbar_3,\Typpialbar_1+\Typpialbar_2+\Typpialbar_3,
\Typpialbar_1+2\Typpialbar_2+\Typpialbar_3\}$.
\newline\newline
$(\Typpibar 4)$ (Extra-rank-4-type) Let $\Typq\in\bKtinf$. Assume $\Typl=4$.
Let $\TypfkAoriginalEX:=\{0\}$.
Define $(\,|\,)$ in the same way as that for ${\mathfrak{sl}}(3|2)$-type.
Let $g_1:=g_2:=0$ and $g_3:=g_4:=1$.
Let $\Typchi$ be such that $\Typq_{ii}:=(-1)^{g_i}\Typq^{(\Typpial_i|\Typpial_i)}$
$(i\in\TypfkI)$ and 
$\Typq_{i^\prime,j^\prime}\Typq_{j^\prime,i^\prime}:=(-1)^{g_{i^\prime}g_{j^\prime}}
\Typq^{2(\Typpial_{i^\prime}|\Typpial_{j^\prime})}$
$(i^\prime,j^\prime\in\TypfkI,\, i^\prime\ne j^\prime)$.
Let $\Typvarpi_i$ be the same as that for ${\mathfrak{sl}}(3|2)$-type.
We have $\TypRchibarpipre=\{\Typpialbar_1, \Typpialbar_1+\Typpialbar_2,
\Typpialbar_2, \Typpialbar_4, \Typpialbar_1+2\Typpialbar_2+3\Typpialbar_3+2\Typpialbar_4,
\Typpialbar_1+2\Typpialbar_2+3\Typpialbar_3+\Typpialbar_4\}$
and $\TypRchibarpipnu=\{\beta\in\TypRchibarpip|\Typchi(\beta,\beta)=-1\}
=\{\Typpialbar_1+\Typpialbar_2+\Typpialbar_3+\Typpialbar_4,
\Typpialbar_1+\Typpialbar_2+\Typpialbar_3,
\Typpialbar_1+2\Typpialbar_2+2\Typpialbar_3+\Typpialbar_4,
\Typpialbar_1+\Typpialbar_2+2\Typpialbar_3+\Typpialbar_4,
\Typpialbar_2+\Typpialbar_3+\Typpialbar_4,
\Typpialbar_3+\Typpialbar_4,
\Typpialbar_2+2\Typpialbar_3+\Typpialbar_4,
\Typpialbar_2+\Typpialbar_3,
\Typpialbar_3\}$.
\newline\newline
$(\Typpibar 5)$ (Extra-rank-2-type) Let $\Typq\in\bKtinf$. Assume $\Typl=2$.
Let $\TypfkAoriginalEX:=\{0\}$.
Define $(\,|\,)$ in the same way as that for ${\mathfrak{sl}}(2|1)$-type.
Let $g^\prime_1:=0$ and $g^\prime_2:=1$.
Let $\zeta\in\bKt$ be such that $\zeta^2+\zeta+1=0$.
Let $\Typchi$ be such that $\Typq_{ii}:=\zeta^{g^\prime_i}\Typq^{(\Typpial_i|\Typpial_i)}$
$(i\in\TypfkI)$ and 
$\Typq_{i^\prime,j^\prime}\Typq_{j^\prime,i^\prime}:=\zeta^{2g^\prime_{i^\prime}g^\prime_{j^\prime}}
\Typq^{2(\Typpial_{i^\prime}|\Typpial_{j^\prime})}$
$(i^\prime,j^\prime\in\TypfkI)$.
Then $\Typvarpi_i$ can be the same as that for ${\mathfrak{sl}}(2|1)$-type.
We have $\TypRchibarpipre=\{\Typpialbar_1, \Typpialbar_1+2\Typpialbar_2\}$
and $\TypRchibarpipnu=\{\Typpialbar_1+\Typpialbar_2,\Typpialbar_2\}$.

\begin{remark}\label{remark:relofInt}
Let $\Typpibar$ be of ${\mathfrak{b}}$-type in $(\Typpibar 2)$,
where ${\mathfrak{b}}$ is one of the finite-dimensional Lie superalgebras of type $A$-$G$
of rank $=\Typl$, see Introduction.
Then ${\mathfrak{b}}$ is identified with ${\mathfrak{g}}$
of Introduction for
$x_{ij}=(\al_i,\al_j)$ and $p(\al_i)=0$ $(b_{ii}\ne -1)$,
$p(\al_j)=1$ $(b_{ii}=-1)$.
Under this identification
(more precisely, under the identification $\{\Typpibar(i)|i\in\TypfkI\}$ and 
$\Typpialbar_i=\Typpibar(i)$ with $\Pi$ and $\al_i$),
$U_{{\acute{Q}}}$ and $\TypIntKhR^\TypIntacuteQ_+$ of Introduction is identified with
$\TypU(\Typchi_{|\TypfkAoriginalpi\times\TypfkAoriginalpi},\Typpibar)$
and $\TypRchibarpip$ respectively
for $\bK=\bC$,
and $\TypRchibarpip$, $\TypRchibarpipre$ and $\TypRchibarpipnu$ are identified with 
$\TypIntRpreal\cup\TypIntRpnull$, $\TypIntRpreal$ and $\TypIntRpnull$ of Introduction
respectively,
see \cite[Lemma~5.6]{AYY15}.
For $\Typpi$ of \eqref{eqn:assumption-a}, we have
$\TypRchipi=\TypRchibarpip\cup(-\TypRchibarpip)$,
$\TypRchipire=\TypRchibarpipre\cup(-\TypRchibarpipre)$,
$\TypRchipinu=\TypRchibarpipnu\cup(-\TypRchibarpipnu)$,
and 
$\TypRchipip=\TypRchipi\cap(\oplus_{i\in\TypfkI}\bZgeqo\Typpial_i)$,
$\TypRchipipre=\TypRchipire\cap(\oplus_{i\in\TypfkI}\bZgeqo\Typpial_i)$,
$\TypRchipipnu=\TypRchipinu\cap(\oplus_{i\in\TypfkI}\bZgeqo\Typpial_i)$,
where recall $\Typpial_i=\Typpi(i)$.

\end{remark}

We see:
\begin{lemma}\label{lemma:asschK}
$\forall \lambda\in\TypfkAoriginal(=\TypfkAoriginalpi
\oplus\TypfkAoriginalEX)$, $\forall\beta\in\TypRchipire$, $\exists k\in\bZ$, $\Typchi(\lambda,\beta)\Typchi(\beta,\lambda)=\Typqbeta^k$. 
\end{lemma}
{\it{Proof.}}
Let $\al^\prime_i:=\Typtau^\Typchi_{f,t}(\Typpibar)(i)$ 
$(i\in\TypfkI)$
and $q^\prime_{i,j}:=\Typchi(\al^\prime_i,\al^\prime_j)$ 
$(i^\prime,j^\prime\in\TypfkI)$.
By \cite[Tables~1-4]{Hec09}, we see that
for $i, j\in\TypfkI$ with $q^\prime_{i,i}\in\bKtinf$,
we have
$q^\prime_{i,j}q^\prime_{j,i}=(q^\prime_{i,i})^{-k}$
for some $k\in\bZgeqo$.
Then the claim for $\lambda\in\TypfkAoriginalpi$ follows form \eqref{eqn:elsRp}.
If $\Typlp=\Typl+1$, The claim for $\lambda\in\TypfkAoriginalEX$
can be seen in a direct way.
\hfill $\Box$
\newline\par
Let $\beta\in\TypRchipire$.
By Lemma~\ref{lemma:asschK},
we have the $\bZ$-module automorphism 
$\Typsbeta:\TypfkAoriginal\to\TypfkAoriginal$
defined in the way that $\Typsbeta(\lambda):=\lambda-k\beta$
for $\lambda\in\TypfkAoriginal$ with $k\in\bZ$ satisfying
$\Typchi(\beta,\lambda)\Typchi(\lambda,\beta)=\Typq_\beta^k$.
We directly see that
$\Typchi(\Typsbeta(\mu),\Typsbeta(\mu))=\Typchi(\mu,\mu)$
and $\Typchi(\Typsbeta(\mu),\Typsbeta(\nu))\Typchi(\Typsbeta(\nu),\Typsbeta(\mu))=\Typchi(\mu,\nu)\Typchi(\nu,\mu)$
for all $\mu$, $\nu\in\TypfkAoriginal$. Note $\Typs_{-\beta}=\Typs_\beta$
and $\Typs_\beta^{-1}=\Typs_\beta$.
If $\beta=\Typpi(i)$ for some $i\in\TypfkI$,
we have $(\Typsbeta)_{|\TypfkAoriginalpi}=\Typacsi$.

Let $\TypWchipire$ be the subgroup of ${\mathrm{Aut}}_\bZ(\TypfkAoriginal)$
generated by 
$\{\Typsbeta|\beta\in\TypRchipire\}$.

\begin{lemma}\label{lemma:asschKb}
{\rm{(1)}} The map 
$\varphi:\TypfkAoriginalpi\to(\bKt)^\Typl$
defined by 
\begin{equation*}
\varphi(\lambda):=(\Typchi(\Typvarpi_j,\lambda)\Typchi(\lambda,\Typvarpi_j))_{j\in\TypfkI}
\quad (\lambda\in\TypfkAoriginalpi)
\end{equation*}
is injective. \newline
{\rm{(2)}} As for
$(\,|\,)$,
the following {\rm{(2a)}} and {\rm{(2b)}}
hold. \par
{\rm{(2a)}} $(\beta|\beta)\ne 0$ and 
$\Typsbeta(\lambda)=\lambda-{\frac {2(\beta|\lambda)} {(\beta|\beta)}}\beta$
for all $\beta\in\TypRchipire$ and all $\lambda\in\TypfkAoriginal$. \par
{\rm{(2b)}} Let $Y_\pm:=\{\beta\in\TypRchipire|\pm(\beta|\beta)>0\}$ and $Z_\pm:=\TyprmSpan_\bR Y_\pm$.
Then $Z_+\cap Z_-=\{0\}$ and $(Z_+|Z_-)=\{0\}$.
Moreover $(Y_+,(\,|\,)_{|Z_+\times Z_+})$ and $(Y_-,-(\,|\,)_{|Z_-\times Z_-})$
are root systems in the sense of  {\rm{\cite[Subsecttion~9.2]{Hum}}}.
\end{lemma} 
{\it{Proof.}}
Recalling \eqref{eqn:pibarpi}, we can directly prove
the claims.
\hfill $\Box$
\newline\par
We also regard $\TypRchipire$ as the root system
in the sense of \cite[Subsection~9.2]{Hum}
with the inner product $(\,|\,)^\prime$
defined by
$(z_1^\pm|z_2^\pm)^\prime=\pm(z_1^\pm|z_2^\pm)$, 
$(z_1^\pm|z_2^\mp)=0$
$(z_1^\pm,z_2^\pm\in Z_\pm)$.
\begin{lemma}\label{lemma:colfcW}
Let $f:\TypfkAoriginalbR\to\bR$ be an $\bR$-linear homomorphism
such that $0\notin f(\TypRchipipre)$. 
Let $X_f^\prime:=\{\beta\in\TypRchipire|f(\beta)>0\}$.
Let $X_f^{\prime\prime}:=\{\beta\in X^\prime|\exists k\in\TypfkJ_{2,\infty}, \gamma_t\in X^\prime
(t\in\TypfkJ_{1,k}), \beta=\sum_{t=1}^k\gamma_t\}$.
Let $X_f:=X_f^\prime\setminus X_f^{\prime\prime}$.
{\rm{(}}Note that $X_f^\prime$ for $f$ with
$f(\Typpial_i)=f(\Typepsilon_r)=1$ equals $\TypRchipipre$.{\rm{)}}
Then we have the following. \newline
{\rm{(1)}} 
$X_f$ is a base of the root system $\TypRchipire$. \newline
{\rm{(2)}} $\TypWchipire$ is the Weyl group of the root system $\TypRchipire$.
As a group, 
$\TypWchipire$ can be presented by the
generators $\{\Typsbeta|\beta\in X_f\}$
and the relations $(\Typsal\Typsbeta)^{m(\al,\beta)}=e$
$(\al,\beta\in X)$, where $m(\al,\beta):=|X_f^\prime\cap(\bR\al+\bR\beta)|$. \newline
{\rm{(3)}} Let $n\in\bN$ and $\gamma_t\in X_f$ $(t\in\TypfkJ_{1,n})$.
Let $w=\Typs_{\gamma_1}\Typs_{\gamma_2}\cdots\Typs_{\gamma_n}$.
Assume that $w\notin\{\Typs_{\beta_1}\Typs_{\beta_2}\cdots\Typs_{\beta_m}|m\in\TypfkJ_{1,n-1},
\beta_r\in X_f (r\in\TypfkJ_{1,m})\}\cup\{e\}$.
Then $\Typs_{\gamma_1}\Typs_{\gamma_2}\cdots\Typs_{\gamma_{t-1}}(\gamma_t)\in X_f^\prime$
for $t\in\TypfkJ_{1,n}$. \newline
{\rm{(4)}} 
Let $\Typhatdelta:=\sum_{\al\in X_f^\prime}\al(\in\TypfkAoriginalpi)$. Then
we have $\Typsal(\Typhatdelta)=\Typhatdelta-2\al$
for all $\al\in X_f$.
Moreover the map $g:\TypWchipire\to\TypfkAoriginalpi$
defined by $g(w):=w(\Typhatdelta)$ is injective.\newline
{\rm{(5)}} We have the group monomorphism
$h:\TypWchipire\to{\mathrm{Aut}}_\bZ(\TypfkAoriginalpi)$
defined by $h(w):=w_{|\TypfkAoriginalpi}$
$(w\in\TypWchipire)$.
\end{lemma}
{\it{Proof.}}
By Lemma~\ref{lemma:asschKb}~(2),
we see that the claims follow from the well-known facts
of the root systems and the Weyl groups.
\hfill $\Box$
\newline\par
By Lemmas~\ref{lemma:colfcW}~(2) and \ref{lemma:asschKb}~(2a), we see that
there exists a unique group homomorphism
$\Typsgnchipi:\TypWchipire\to\{-1,1\}$
such that $\Typsgnchipi(\Typsbeta)=-1$
$(\beta\in\TypRchipipre)$.

\subsection{Irreducible modules}
Let $\TypLam:\TypUochipi\to\bK$ be a $\bK$-algebra homomorphism.
Then there exists an irreducible
$\TypUchipi$-module $\TypmclLchipiLam$
such that there exists $v_\TypLam\in\TypmclLchipiLam\setminus\{0\}$
satisfying the conditions that letting $\TypmclLchipiLam_\lambda:=\TypUchipi_\lambda v_\TypLam$
$(\lambda\in\TypfkAoriginalpi)$, we have
\begin{equation}\label{eqn:codLLam}
\begin{array}{l}
\TypmclLchipiLam=\oplus_{\nu\in\TypfkAoriginalpip}\TypmclLchipiLam_{-\nu},\,\,
\TypmclLchipiLam_\mu=\{0\}\,(-\mu\notin\TypfkAoriginalpip),\\
\TypmclLchipiLam_{-\nu}=\TypUmchipi_{-\nu}v_\TypLam
\,(\nu\in\TypfkAoriginalpip),
\\
\mbox{$Zv_\TypLam=\TypLam(Z)v_\TypLam$
$(Z\in\TypUochipi)$
and $\TyptrE_iv_\TypLam=0$ ($i\in\TypfkI$).}
\end{array}
\end{equation}
There also exists a 
$\TypUchipi$-module $\TypmclMchipiLam$
such that there exists ${\tilde v}_\TypLam\in\TypmclMchipiLam\setminus\{0\}$
satisfying the same conditions of those of \eqref{eqn:codLLam}
with $\TypmclMchipiLam$ and ${\tilde v}_\TypLam$ in place of 
$\TypmclMchipiLam$ and $v_\TypLam$ respectively
and the condition that $\dim\TypUmchipi_{-\nu}{\tilde v}_\TypLam=\dim\TypUmchipi_{-\nu}$
$(\nu\in\TypfkAoriginalpip)$.
Then there exists a unique $\TypUchipi$-module epimorphism
$f:\TypmclMchipiLam\to\TypmclLchipiLam$ with
$f({\tilde v}_\TypLam)=v_\TypLam$.
It is clear that $f(\TypmclMchipiLam_\lambda)=\TypmclLchipiLam_\lambda$
$(\lambda\in\TypfkAoriginalpi)$.
We call $\TypmclLchipiLam$ (resp. $\TypmclMchipiLam$)
{\it{the highest weight irreducible {\rm{(}}resp. the Verma{\rm{)}} $\TypUchipi$-module}}
of $\TypLam$.
For $i\in\TypfkI$, we can easily see that
\begin{equation}\label{eqn:Fimv}
\TyptrF_i^mv_\TypLam\ne 0\quad\Longleftrightarrow\quad
(m)_{q_{ii}}!(m;q_{ii}^{-1};\TypLam(\TyptrK_{\Typpial_i}\TyptrL_{-\Typpial_i}))!\ne 0.
\end{equation}

Define the map $\TypOchipi:\TypUochipi\to\bK$ by
$\TypOchipi(Z):=0$ for all $(Z\in\TypUochipi)$.

Let $i\in\TypfkI$. Define the map $\TyptauichipiLam:\TypUochipi\to\bK$ as follows.
If there exists $h\in\bZgeqo$ such that $\TyptrF_i^hv_\TypLam\ne 0$
and $\TyptrF_i^{h+1}v_\TypLam=0$,
let $\TyptauichipiLam:\TypUochipi\to\bK$ be
the $\bK$-algebra homomorphism defined
by $\TyptauichipiLam(\TyptrK_\lambda\TyptrL_\mu)=
{\frac {\Typchi(h\Typpial_i,\mu)} {\Typchi(\lambda,h\Typpial_i)}}
\TypLam(\TyptrK_\lambda\TyptrL_\mu)$
$(\lambda,\mu\in\TypfkAoriginal)$.
If $\TyptrF_i^mv_\TypLam\ne 0$ for all $m\in\bN$,
let $\TyptauichipiLam:=\TypOchitauipi$.
We also let $\Typtauichipi(\TypOchipi):=\TypOchitauipi$.

If $\TyptauichipiLam\ne\TypOchitauipi$, letting
$\Typpi^\prime:=\Typchitauipi$ and $\TypLam^\prime:=\TyptauichipiLam$, 
we have the $\bK$-linear isomorphism
\begin{equation*}
\mbox{$\TyphatTchiPRpiPRLami:\TypmcLchiPRpiPRLami\to\TypmclLchipiLam$
with
$\TyphatTchiPRpiPRLami(Xv_{\TypLam^\prime})=\TypTchiPRpii(X)\TyptrF_i^hv_\TypLam$
$(X\in\TypUchiPRpi)$,}
\end{equation*} 
where let $h$ be as above, see \cite[Lemma~6.3]{AYY15}. 

For $f\in\TypbBI$ and $t\in\bZgeqo$, define 
$\TypchitauftTyppiLam$ in the way that 
$\TypchitaufoTyppiLam:=\TypLam$ 
and $\TypchitauftTyppiLam:=\Typtau^{\Typchi,\Typpi^\prime}_{f(t)}(\TypchitauftmTyppiLam)$
(if $t\in\bN$),
where $\Typpi^\prime:=\TypchitauftmTyppi$.

It follows from \cite[Lemma~6.6]{AYY15} that
\begin{equation*}
\mbox{$\dim\TypmclLchipiLam<\infty\quad\Longleftrightarrow
\quad\forall f\in\TypbBI,\,\forall t\in\bZgeqo, \TypchitauftTyppiLam\ne\TypOchitauipift$.}
\end{equation*} 
\subsection{Main theorem for irreducible cases}
Let $\Typc_i:=\Typkpch(\Typq_{ii})\,(=\Typc_{\Typpial_i})$
$(i\in\TypfkI)$. Then we have
\begin{equation}\label{eqn:rcrho}
\Typchi(\Typpial_i,\lambda)^{\Typc_i-1}\Typchi(\lambda,\Typpial_i)^{\Typc_i-1}=
{\frac{\Typhrhochitauipi(\lambda)} {\Typhrhochipi(\lambda)}}
\quad(i\in\TypfkI,\,\lambda\in\TypfkAoriginalpi),
\end{equation} see \cite[Lemma~7.2]{BY18} for example.

For $\beta\in\TypRchipipre$, let
$\TypRchipibeta:=\{\al\in\TypRchipip|\Typsbeta(\al)\in -\TypRchipip\}$.
\begin{proposition}\label{proposition:mqbeta}
Let $\beta\in\TypRchipipre$.
Then 
\begin{equation}\label{eqn:eqmqb}
\mbox{$\exists m\in\bZ$, $\Typhrhochipi(\beta)=\Typqbeta^m$\,\, and\,\,
$m\beta=\sum_{\Typpial\in\TypRchipibeta}(1-\Typcal)\al$.}
\end{equation} 
\end{proposition}
{\it{Proof.}}
Let $i\in\TypfkI$. Recall $\Typpial_i=\Typpi(i)$.
Recall from \eqref{eqn:tautaupi} that $\TypRchipiptaui=\TypRchipip\setminus\{\Typpial_i\}\cup\{-\Typpial_i\}$.
If $\beta=\Typpial_i$, we have \eqref{eqn:eqmqb} for $m=1$
since $\Typc_i=0$.
We show \eqref{eqn:eqmqb} for $\Typchitauipi$ in place of $\Typpi$
after assuming $\eqref{eqn:eqmqb}$ for $\Typpi$ is correct.
Recall that for $\beta\in\TypRchipip$, we have
$\Typsbeta(\Typpial_i)=\Typpial_i-n\beta$
with $n\in\bZ$ satisfying $\Typqbeta^n=
\Typchi(\beta,\Typpial_i)\Typchi(\Typpial_i,\beta)$.

For $\beta\in\TypRchipiptaui$, 
we see that
\begin{equation*}
\TypRchitauipibeta=
\left\{\begin{array}{ll}
\{-\Typpial_i\} & \quad\mbox{if $\beta=-\Typpial_i(=\Typchitauipi(i)))$,} \\
\TypRchipibeta & \quad\mbox{if $n=0$,} \\
\TypRchipibeta\cup\{\Typsbeta(\Typpial_i),-\Typpial_i\} & \quad\mbox{if $n\ne 0$ and $\Typsbeta(\Typpial_i)\in\TypRchipip$,} \\
\TypRchipibeta\setminus\{-\Typsbeta(\Typpial_i),\Typpial_i\} & \quad\mbox{if $n\ne 0$, 
$\beta\ne-\Typpial_i$ and $\Typsbeta(\Typpial_i)\in-\TypRchipip$}
\end{array}\right.
\end{equation*}
Then we easily obtain \eqref{eqn:eqmqb}
from \eqref{eqn:rcrho}.
\hfill $\Box$
\newline\par
Define the $\bZ$-module isomorphism
$\Typhrhochipifull:\TypfkAoriginal\to\bKt$ by
$(\Typhrhochipifull)_{|\TypfkAoriginalpi}:=\Typhrhochipi$
and $\Typhrhochipifull(\TypfkAoriginalEX)=\{1\}$.

Define the subalgebra $\TypUodagchipi$ of $\TypUochipi$
by $\TypUodagchipi:=\oplus_{\lambda\in\TypfkAoriginal}\bK
\TyptrK_\lambda\TyptrL_{-\lambda}$.
Define the $\bK$-algebra automorphism
$\Typhrhodag:\TypUodagchipi\to\TypUodagchipi$ by
$\Typhrhodag(\TyptrK_\lambda\TyptrL_{-\lambda})
:=\Typhrhochipifull(\lambda)\TyptrK_\lambda\TyptrL_{-\lambda}$
$(\lambda\in\TypfkAoriginal)$.
For $\beta\in\TypRchipipre$, define the $\bK$-algebra automorphism
$\Typsdagbeta:\TypUodagchipi\to\TypUodagchipi$ by
$\Typsdagbeta(\TyptrK_\lambda\TyptrL_{-\lambda}):=\TyptrK_{\Typs_\beta(\lambda)}
\TyptrL_{-\Typs_\beta(\lambda)}$
$(\lambda\in\TypfkAoriginal)$.
Let
$\Typdotsdagbeta:=\Typhrhodag\circ\Typsdagbeta\circ\Typhrhodag^{-1}$.
Then $\Typdotsdagbeta(\TyptrK_\lambda\TyptrL_{-\lambda})=\Typhrhochipibeta^{-t}\TyptrK_{\Typs_\beta(\lambda)}\TyptrL_{-\Typs_\beta(\lambda)}$
for $\lambda\in\TypfkAoriginal$
and $t\in\bZ$ with
$\Typchi(\beta,\lambda)\Typchi(\lambda,\beta)=\Typq_\beta^t$.
Let $\TypWchipidagre$ (resp. $\TypdotWchipidagre$) be the subgroup of ${\mathrm{Aut}}_{\bK-{\mathrm{algebra}}}(\TypUodagchipi)$ generated by $\{\Typsdagbeta|\beta\in\TypRchipipre\}$
(resp. $\{\Typdotsdagbeta|\beta\in\TypRchipipre\}$). Obviously we have 
\begin{equation*}
\begin{array}{l}
\mbox{the group isomorphism $\Typxichipidag:\TypWchipire\to\TypWchipidagre$
(resp. $\Typdotxichipidag:\TypWchipire\to\TypdotWchipidagre$)} \\
\mbox{defined by $\Typxichipidag(\Typsbeta)=\Typsdagbeta$
(resp. $\Typdotxichipidag(\Typsbeta)=\Typdotsdagbeta$)
for $\beta\in\TypRchipipre$.}
\end{array}
\end{equation*} 

Let $\Typomegazero:\TypfkAoriginal\to\bKt$ be the $\bK$-algebra isomorphism defined by
$\Typomegazero(\TypfkAoriginal)=\{1\}$.
Let $\TypfkBchipionedagpre:=\{Z\in\TypUodagchipi|\forall\beta\in\TypRchipipre,
\Typdotsdagbeta(Z)=Z\}(=\{Z\in\TypUodagchipi||\forall w\in\TypdotWchipidagre, w(Z)=Z\})$.
By Lemma~\ref{lemma:asschK}, we see that for $Z\in\TypUodagchipi$,
\begin{equation*}
Z\in\TypfkBchipionedagpre\quad\Longleftrightarrow\quad
\forall\beta\in\TypRchipipre,\,\mbox{$Z$ satisfies $(e1)_\beta$ and $(e2)_\beta$.}
\end{equation*}

For $\beta\in\TypRchipipnu$, we
let 
\begin{equation*}
\TypPchipipmbeta:=\sum_{t=0}^{\Typcbeta-1}(\Typhrhochipibeta\TyptrK_\beta\TyptrL_{-\beta})^{\pm t}
\in\TypUodagchipi,
\end{equation*} and we have $\TypPchipimbeta=(\Typhrhochipibeta\TyptrK_\beta\TyptrL_{-\beta})^{1-\Typcbeta}\TypPchipipbeta$,
where if $\Typqbeta=1$, $\TypPchipipmbeta=0$ since $\Typcbeta=0$.
Let $\TyphatPchipi:=\prod_{\beta\in\TypRchipipnu}\TypPchipipbeta\TypPchipimbeta$.
We have  $\TyphatPchipi\in\TypfkBchipionedagpre$.
Moreover we can easily see the following.
\begin{lemma} \label{lemma:KEY} We have
\begin{equation*}
\TyphatPchipi\cdot\TypfkBchipionedagpre\subset\TypUodagchipi\cap\TypfkBchipiomegazero.
\end{equation*}
\end{lemma}

\begin{equation}\label{eqn:assTyp}
\begin{array}{l}
\mbox{From now on until Remark~\ref{remark:rankonecase},} \\
\mbox{we assume that $\dim\TypmclLchipiLam<\infty$
and $\TypLam(\TyphatPchipi)\ne 0$.}
\end{array}
\end{equation}
Then we call $\TypmclLchipiLam$ {\it{typical}}.
Let $\TypmclELam$ be a free $\bZ$-module with a basis $\{\TypeLam_\nu|\nu\in\TypfkAoriginalpi\}$,
i.e., $\TypmclELam=\oplus_{\nu\in\TypfkAoriginalpi}\bZ\TypeLam_\nu$.
For $\nu\in\TypfkAoriginalpi$, define the $\bK$-algebra homomorphism
$\TypLamnu:\TypUochipi\to\bK$ by 
$\TypLamnu(\TyptrK_\lambda\TyptrL_\mu):=
{\frac {\Typchi(\lambda,\nu)} {\Typchi(\nu,\mu)}}\TypLam(\TyptrK_\lambda\TyptrL_\mu)$
$(\lambda,\mu\in\TypfkAoriginal)$. 
%By \cite{Hec09}, and by using Weyl groupoids \cite{HY08}, \cite{Y16},
%we see that
%\begin{equation*}
%\forall\beta\in\TypRchipipre, \forall\nu\in\TypfkAoriginalpi,
%\exists k\in\bZ, \Typchi(\beta,\nu)\Typchi(\lambda,\nu)=\Typqbeta^k.
%\end{equation*} 

By Lemma~\ref{lemma:asschKb}~(1), we see that for $\mu$, $\nu\in\TypfkAoriginalpi$, 
\begin{equation}\label{eqn:ppreinv}
\mu=\nu\,\,\Leftrightarrow\,\,\TypeLam_\mu=\TypeLam_\nu\,\,\Leftrightarrow\,\,
(\TypLammu)_{|\TypUodagchipi}=(\TypLamnu)_{|\TypUodagchipi}\,\,\Leftrightarrow\,\,
\TypLammu=\TypLamnu.
\end{equation} 

Since $\dim\TypmclLchipiLam<\infty$,
by \eqref{eqn:Fimv} and \eqref{eqn:eqmqb},
we have 
\begin{lemma}
Let $\beta\in\TypRchipipre$ and $\nu\in\TypfkAoriginalpi$.
Then there exist $a_\beta^\Lambda(\nu)$, $b_\beta^\Lambda(\nu)\in\bZ$ such that
\begin{equation}\label{eqn:defab}
\TypLamnu(\TyptrK_\beta\TyptrL_{-\beta})=\Typqbeta^{a_\beta^\Lambda(\nu)}
={\frac {\Typqbeta^{b_\beta^\Lambda(\nu)}} {\Typhrhochipibeta}}.
\end{equation} 
\end{lemma}

Let $\TypBifkApi$ be the set of bijections from $\TypfkAoriginalpi$
to $\TypfkAoriginalpi$.

Let $\beta\in\TypRchipipre$. Define
$\TypsigLam_\beta\in\TypmBiclELam$
(resp. $\TypdotsigLam_\beta\in\TypmBiclELam$)
in the way that $\TypsigLam_\beta(\TypeLam_\nu):=\TypeLam_{\nu-a_\beta^\Lambda(\nu)\beta}$
(resp. $\TypdotsigLam_\beta(\TypeLam_\nu):=\TypeLam_{\nu-b_\beta^\Lambda(\nu)\beta}$)
for $\nu\in\TypfkAoriginalpi$,
Let $r_\beta\in\bZ$ be such that 
\begin{equation}\label{eqn:rbeta}
\Typqbeta^{r_\beta}=\Typhrhochipibeta,
\end{equation} see \eqref{eqn:eqmqb}.
Then we see that $b_\beta^\Lambda(\nu)=a_\beta^\Lambda(\nu)+r_\beta$ and
\begin{equation}\label{eqn:sbnemu}
\TypsigLam_\beta(\TypeLam_{\nu+\mu})=\TypeLam_{\nu-a_\beta^\Lambda(\nu)\beta+\Typsbeta(\mu)}
\quad\TypdotsigLam_\beta(\TypeLam_{\nu+\mu})=\TypeLam_{\nu-(a_\beta^\Lambda(\nu)+r_\beta)\beta+\Typsbeta(\mu)}\quad 
(\mu, \nu\in\TypfkAoriginalpi).
\end{equation}
Define $\TypsLam_\beta\in\TypBifkApi$
(resp. $\TypdotsLam_\beta\in\TypBifkApi$) by
$\TypsLam_\beta(\mu):=\Typsbeta(\mu)-a_\beta^\Lambda(0)\beta$
(resp. 
$\TypdotsLam_\beta(\mu):=\Typsbeta(\mu)-(a_\beta^\Lambda(0)+r_\beta)\beta$
for $\mu\in\TypfkAoriginalpi$.
We have $\TypsigLam_\beta(\TypeLam_\mu)=\TypeLam_{\TypsLam_\beta(\mu)}$
and $\TypdotsigLam_\beta(\TypeLam_\mu)=\TypeLam_{\TypdotsLam_\beta(\mu)}$
for $\mu\in\TypfkAoriginalpi$.

We have
\begin{equation}\label{eqn:preinv}
\begin{array}{l}
\forall\nu\in\TypfkAoriginalpi, \forall\beta\in\TypRchipipre, \forall Z\in\TypUodagchipi, \\
\quad\TypLam^{+\TypsLam_\beta(\nu)}(Z)=\TypLamnu(\Typsdagbeta(Z)),
\,\,\TypLam^{+\TypdotsLam_\beta(\nu)}(Z)=\TypLamnu(\Typdotsdagbeta(Z)).
\end{array}
\end{equation} 

Let $\TypfkWchipireLam$ (resp. $\TypdotfkWchipireLam$) 
be the subgroup of $\TypmBiclELam$
generated by $\{\TypsigLam_\beta|\beta\in\TypRchipipre\}$
(resp. $\{\TypdotsigLam_\beta|\beta\in\TypRchipipre\}$).
Let $\TypWchipireLam$ (resp. $\TypdotWchipireLam$) 
be the subgroup of $\TypBifkApi$
generated by $\{\TypsLam_\beta|\beta\in\TypRchipipre\}$
(resp. $\{\TypdotsLam_\beta|\beta\in\TypRchipipre\}$).
We have the group isomorphism
$\TypetachipiLamp:\TypfkWchipireLam\to\TypWchipireLam$
(resp. $\TypdotetachipiLamp:\TypdotfkWchipireLam\to\TypdotWchipireLam$)
with $\TypetachipiLamp(\TypsigLam_\beta)=\TypsLam_\beta$
(resp. $\TypdotetachipiLamp(\TypsigLam_\beta)=\TypdotsLam_\beta$)
for $\beta\in\TypRchipipre$.

\begin{lemma}
We have
the group isomorphism $\TypetachipiLam:\TypWchipire\to\TypWchipireLam$
{\rm{(}}resp. $\TypdotetachipiLam:\TypWchipire\to\TypdotWchipireLam${\rm{)}}
defined by $\TypetachipiLam(\Typsbeta)=\TypsLam_\beta$
{\rm{(}}resp. $\TypdotetachipiLam(\Typsbeta)=\TypdotsLam_\beta${\rm{)}}
for $\beta\in\TypRchipipre$.
\end{lemma}
{\it{Proof.}} We prove the claim for $\TypetachipiLam$.
We can similarly prove the claim for $\TypdotetachipiLam$.
By \eqref{eqn:ppreinv} and \eqref{eqn:preinv}, the group epimorphism 
$\TypetachipiLam$ exists, where we use $\Typdotxichipidag$.
By Lemma~\ref{lemma:colfcW}~(4) and \eqref{eqn:sbnemu}, we see that
for sufficiently large $k\in\bN$,
the map $g:\TypWchipire\to\TypfkAoriginalpi$
defined by $g(w):=\TypetachipiLam(w)(k\Typhatdelta)$ is injective.
\hfill $\Box$
\newline\par
We have the group isomorphisma $\TypxichipiLam:=(\TypetachipiLamp)^{-1}\circ\TypetachipiLam:\TypWchipire
\to\TypfkWchipireLam$ and
$\TypdotxichipiLam:=(\TypdotetachipiLamp)^{-1}\circ\TypdotetachipiLam:\TypWchipire\to\TypdotfkWchipireLam$.

As for the following proposition,
compare with \cite[Lemma and its proof in Subsection~23.3]{Hum}.

\begin{proposition}\label{proposition:Znuprime}
 Let $\nu\in\TypfkAoriginalpi$. Then it follows that
\begin{equation}\label{eqn:trueinv}
\forall Z\in\TypUodagchipi\cap\TypfkBchipiomegazero, \TypLamnu(Z)=\TypLam(Z)
\quad\Longleftrightarrow\quad\exists\gamma\in\TypdotWchipireLam, \gamma(\nu)=0.
\end{equation}
\end{proposition}
{\it{Proof.}} The fact in \eqref{eqn:preinv} implies $\Longleftarrow$.

We prove $\Longrightarrow$.
Note $\TypLam=\TypLamo$.
For $\lambda\in\TypfkAoriginalpi$, let 
$X_\lambda:=\{\gamma(\lambda)\in\TypfkAoriginalpi|\gamma\in\TypdotWchipireLam\}$
and ${\mathcal{F}}_\lambda:=\{\TyptrK_{\Typvarpi_j}\TyptrL_{-\Typvarpi_j}-\TypLamlambda(\TyptrK_{\Typvarpi_j}\TyptrL_{-\Typvarpi_j})|j\in\TypfkI\}$.
Assume $0\notin X_\nu$.
Let ${\mathcal{F}}:=(\cup_{\lambda\in X_0\cup X_\nu}{\mathcal{F}}_\lambda)\setminus{\mathcal{F}}_0$.
Let $Z^{\prime\prime}_\nu:=\TyphatPchipi\cdot\prod_{h\in{\mathcal{F}}}h$.
Then $\TypLam(Z^{\prime\prime}_\nu)\ne 0$ and $\TypLamlambda(Z^{\prime\prime}_\nu)=0$
for all $\lambda\in (X_0\cup X_\nu)\setminus\{0\}$. 
Let $Z_\nu^\prime:={\frac 1 {\TypLam(Z^{\prime\prime}_\nu)}}Z^{\prime\prime}_\nu$.
Let $Y:=\{a(Z_\nu^\prime)|a\in\TypdotWchipidagre\}$,
and $Z_\nu:=\sum_{{\acute{Z}}\in Y}{\acute{Z}}$. 
We have $Z_\nu\in\TypUodagchipi\cap\TypfkBchipiomegazero$ and
$\TypLamnu(Z_\nu)=0$.
Let $y:=|Y|$.
Let $a_t\in\TypdotWchipidagre$ $(t\in\TypfkJ_{1,y})$ be
such that  $Y:=\{a_t(Z_\nu^\prime)|t\in\TypfkJ_{1,y}\}$.
We may asumme $a_1=e$.
Let $\gamma_t:=(\TypdotetachipiLam\circ(\Typdotxichipidag)^{-1})(a_t^{-1})$ $(t\in\TypfkJ_{1,y})$.
Then
$\TypLam(Z_\nu)=|\{t\in\TypfkJ_{1,y}|\gamma_t(0)=0\}|\geq 1$.
\hfill $\Box$

\begin{remark}
Keep the notation in Proof of Proposition~\ref{proposition:Znuprime}.
Let ${\dot G}_\dagger:=\{a\in\TypdotWchipidagre|a(Z_\nu^\prime)=Z_\nu^\prime\}$
and ${\dot G}_\TypLam:=\{\gamma\in\TypdotWchipireLam|\gamma(0)=0\}$.
Let $y^\prime:=|\{t\in\TypfkJ_{1,y}|\gamma_t(0)=0\}|$,
and assume that $\gamma_{t^\prime}(0)=0$ $(t^\prime\in\TypfkJ_{1,y^\prime})$
and $\gamma_{t^{\prime\prime}}(0)\ne 0$ $(t^{\prime\prime}\in\TypfkJ_{{y^\prime}+1,y})$.
Then ${\dot G}_\TypLam=\cup_{t=1}^{y^\prime}
(\TypdotetachipiLam\circ(\TypdotxichipiLam)^{-1})({\dot G}_\dagger a_t^{-1})$,
where the cup of RHS is disjoint union.
In particular, $\TypLam(Z_\nu)=y^\prime={\frac {|{\dot G}_\TypLam|} {|{\dot G}_\dagger|}}$.
\end{remark}

Let $\TyphatmclELam$ be the $\bZ$-module formed by
all $\bZ$-module homomorphisms
from $\TypmclELam$ to $\bZ$.
Then there exists a unique $\bZ$-module 
monomorphism $\TypiotaLam:\TypmclELam\to\TyphatmclELam$
with $\TypiotaLam(\TypeLam_\nu)(\TypeLam_\lambda)=\delta_{\nu,\lambda}$.
We shall denote $\TypiotaLam(\TypeLam_\nu)$ by $\TypeLam_\nu$
for simplicity.
We shall regard $\TypmclELam$ as the $\bZ$-submodule of $\TyphatmclELam$ in a natural way.

Let
\begin{equation*}
[\TypmclLchipiLam]^\TypLam:=\sum_{\nu\in\TypfkAoriginalpip}\dim\TypmclLchipiLam_{-\nu}
\TypeLam_{-\nu}.
\end{equation*}
\begin{lemma}\label{lemma:sbetaLam} We have:
\begin{equation}\label{eqn:sbetaLamfm}
\forall\beta\in\TypRchipipre,\,\,
\TypsigLam_\beta([ \TypmclLchipiLam ]^\TypLam)=[ \TypmclLchipiLam ]^\TypLam.
\end{equation}
\end{lemma}
{\it{Proof.}}
Use the notation of Theorem~\ref{theorem:PBWfinite}~(2).  
Let ${\mathcal{U}}_\beta$ be the $\bK$-algebra of $\TypUchipi$ generated
by $\TyptrK_\beta^{\pm 1}$, $\TyptrL_\beta^{\pm 1}$, 
$\TyptrE_\beta$, $\TyptrF_\beta$.
Let ${\mathbb{X}}:=\bZgeqo\times\bKt$.
For $(n,a)\in{\mathbb{X}}$, 
let $V_{n,a}$ be an irreducible $(n+1)$-dimensional left ${\mathcal{U}}_\beta$-module
such that there exists $v_{n,a}\in V_{n,a}\setminus\{0\}$
with $\TyptrE_\beta v_{n,a}=0$,
$\TyptrK_\beta v_{n,a}=av_{n,a}$
and $\TyptrL_\beta v_{n,a}=a\Typqbeta^{-n}v_{n,a}$.
Since the actions of $\TyptrK_\beta^{\pm 1}$ and $\TyptrL_\beta^{\pm 1}$
on $\TypmclLchipiLam$ are diagonalizable,
the following $(*)$ and $(**)$ hold. 
For a $\bK$-linear subspace $M$ of $\TypmclLchipiLam$
and $x$, $y\in\bKt$
(resp. $\lambda\in\TypfkAoriginalpi$), let $M_{x,y}:=\{m\in M|\TyptrK_\beta m=xm,
\TyptrL_\beta m=ym\}$ (resp. $M_\lambda:=M\cap\TypmclLchipiLam_\lambda$).
Note $\TypmclLchipiLam=\oplus_{x,y\in\bKt}\TypmclLchipiLam_{x,y}$.
\newline\newline
$(*)$ Any ${\mathcal{U}}_\beta$-submodule $M$ of $\TypmclLchipiLam$ satisfies
$M=\oplus_{x,y\in\bKt}M_{x,y}$.
Moreover there exists a $\bK$-linear subspace $N^\prime$ of $M$
such that $M=N^\prime\oplus N$ and 
$N^\prime=\oplus_{x,y\in\bKt}N^\prime_{x,y}$.\newline
$(**)$ For any nonzero ${\mathcal{U}}_\beta$-submodule $M$ of $\TypmclLchipiLam$
and any proper ${\mathcal{U}}_\beta$-submodule $N$ of $M$,
any non-zero irreducible ${\mathcal{U}}_\beta$-submodule
of  $M/N$ is isomorphic to $V_{n,a}$
for some $(n,a)\in{\mathbb{X}}$.\newline
\newline\par
We shall show the following $(\sharp)$.
\newline\newline
$(\sharp)$ As a ${\mathcal{U}}_\beta$-module,
$\TypmclLchipiLam$ is a completely reducible module
whose  irreducible components are isomorphic to $V_{n,a}$ for some $(n,a)\in{\mathbb{X}}$.
\newline\par
Let $Z:=\Typqbeta\TyptrK_\beta-\TyptrL_\beta+(1-\Typqbeta)\TyptrF_\beta\TyptrE_\beta$
and $Y:=\TyptrK_\beta\TyptrL_\beta$. Then 
$Z$ and $Y$ are central elements of ${\mathcal{U}}_\beta$.
Note that $Z v_{n,a}=(a\Typqbeta-a\Typqbeta^{-n})v_{n,a}$
and $Y v_{n,a}=a^2\Typqbeta^{-n}v_{n,a}$.
For $(n,a)$, $(r,b)\in{\mathbb{X}}$, we see that $a\Typqbeta-a\Typqbeta^{-n}=b\Typqbeta-b\Typqbeta^{-r}$
and $a^2\Typqbeta^{-n}=b^2\Typqbeta^{-r}$
if and only if  $(n,a)=(r,b)$.
For $(n,a)\in{\mathbb{X}}$, let ${\mathcal{L}}_{(n,a)}:=\{v\in\TypmclLchipiLam|
Z v=(a\Typqbeta-a\Typqbeta^{-n})v, Yv=a^2\Typqbeta^{-n}v\}$
and  ${\mathcal{L}}^\flat_{(n,a)}:=\{v\in\TypmclLchipiLam|
Yv=a^2\Typqbeta^{-n}v, \exists h\in\bN, (Z-(a\Typqbeta-a\Typqbeta^{-n}))^hv=0\}$.
It is clear that ${\mathcal{L}}_{(n,a)}\subset{\mathcal{L}}^\flat_{(n,a)}$ and 
${\mathcal{L}}^\flat_{(n,a)}\ne \{0\}\Leftrightarrow{\mathcal{L}}_{(n,a)}\ne \{0\}$.
Let ${\mathbb{X}}^\prime:=\{(n,a)\in{\mathbb{X}}|{\mathcal{L}}_{(n,a)}\ne \{0\}\}$.
By $(*)$ and $(**)$, we have $\TypmclLchipiLam=\oplus_{(n,a)\in{\mathbb{X}}^\prime}{\mathcal{L}}^\flat_{(n,a)}$.
Note that ${\mathcal{L}}_{(n,a)}$ and ${\mathcal{L}}^\flat_{(n,a)}$
are ${\mathcal{U}}_\beta$-submodules of $\TypmclLchipiLam$,
and that ${\mathcal{L}}_{(n,a)}=\oplus_{\lambda\in\TypfkAoriginalpi}({\mathcal{L}}_{(n,a)})_\lambda$
and ${\mathcal{L}}^\flat_{(n,a)}=\oplus_{\lambda\in\TypfkAoriginalpi}({\mathcal{L}}^\flat_{(n,a)})_\lambda$.
Note the following.
\newline\newline
$(\sharp\sharp)$ For any proper ${\mathcal{U}}_\beta$-submodule $M$ of ${\mathcal{L}}^\flat_{(n,a)}$,
a non-zero irreducible ${\mathcal{U}}_\beta$-submodule of ${\mathcal{L}}^\flat_{(n,a)}/M$
is isomorphic to $V_{n,a}$.
\newline\par
Let $(n,a)\in{\mathbb{X}}^\prime$
and let ${\mathcal{L}}^{\flat,(t)}_{(n,a)}:=\{v\in{\mathcal{L}}^\flat_{(n,a)}|
\TyptrK_\beta v= a\Typqbeta^{-t}v,
\TyptrL_\beta v= a\Typqbeta^{t-n}v\}
=\{v\in{\mathcal{L}}^\flat_{(n,a)}|
\TyptrK_\beta v= a\Typqbeta^{-t}v,
\TyptrK_\beta\TyptrL_\beta^{-1} v= \Typqbeta^{n-2t}v\}$ $(t\in\bZ)$.
Let ${\mathcal{L}}^{\flat\flat}_{(n,a)}:=\oplus_{t=-\infty}^\infty{\mathcal{L}}^{\flat,(t)}_{(n,a)}$.
Since ${\mathcal{L}}^{\flat\flat}_{(n,a)}$ is a ${\mathcal{U}}_\beta$-submodule of ${\mathcal{L}}^\flat_{(n,a)}$,
we have 
\begin{equation*}
\mbox{${\mathcal{L}}^{\flat\flat}_{(n,a)}=\oplus_{t=0}^n{\mathcal{L}}^{\flat,(t)}_{(n,a)}$
and $\dim{\mathcal{L}}^{\flat,(t)}_{(n,a)}=\dim{\mathcal{L}}^{\flat,(0)}_{(n,a)}$
$(t\in\TypfkJ_{1,n})$.}
\end{equation*}
We have ${\mathcal{L}}^{\flat,(t)}_{(n,a)}=\oplus_{\lambda\in\TypfkAoriginalpi}({\mathcal{L}}^{\flat,(t)}_{(n,a)})_\lambda$.
Let $\mu\in\TypfkAoriginalpi$ be such that $({\mathcal{L}}^\flat_{(n,a)})_\mu\ne\{0\}$.
Let $v\in({\mathcal{L}}^\flat_{(n,a)})_\mu\setminus\{0\}$.
Let $k\in\bZgeqo$
be such that $\TyptrE_\beta^tv\ne 0$ $(t\in\TypfkJ_{0,k})$
and $\TyptrE_\beta^{k+1}v=0$.
We see that $\oplus_{r=0}^n\bK\TyptrF_\beta^r\TyptrE_\beta^kv$
is a ${\mathcal{U}}_\beta$-submodule of ${\mathcal{L}}^\flat_{(n,a)}$
isomorphic to $V_{n,a}$ with $\TyptrE_\beta^kv\in({\mathcal{L}}^{\flat,(0)}_{(n,a)})_{\mu+k\Typpial_i}\setminus\{0\}$.
Hence $v\in({\mathcal{L}}^{\flat,(k)}_{(n,a)})_\mu\setminus\{0\}$.
Hence
we have ${\mathcal{L}}^\flat_{(n,a)}={\mathcal{L}}^{\flat\flat}_{(n,a)}$.
Then we see that \eqref{eqn:sbetaLamfm} is true
(and that $(\sharp)$ is also true). 
\hfill $\Box$
\newline\par
Let $(\TyphatmclELam)^\prime:=\{\sum_{\mu\in\TypfkAoriginalpip}a_\mu\TypeLam_{\nu-\mu}
\in\TyphatmclELam 
|\nu\in\TypfkAoriginalpi,\,\,a_\mu\in\bZ\,(\mu\in\TypfkAoriginalpi)\}$.
Let $\TypcheckmclELam:=\{\sum_{t=1}^kc_t\in\TyphatmclELam|k\in\bN,\,c_t\in (\TyphatmclELam)^\prime\,
(t\in\TypfkJ_{1,k})\}$.
In fact, $\TypcheckmclELam=(\TyphatmclELam)^\prime$.

For $\nu\in\TypfkAoriginalpi$, let
$[\![\nu]\!]^\TypLam:=\sum_{\mu\in\TypfkAoriginalpip}(\dim\TypUmchipi_{-\mu})
\TypeLam_{\nu-\mu}\,(\in\TypcheckmclELam)$, i.e.,
\begin{equation*}
[\![\nu]\!]^\TypLam(\TypeLam_{\nu+\lambda})
=\dim\TypmclM^{\Typchi,\Typpi}(\TypLamnu)_\lambda
\quad (\lambda\in\TypfkAoriginalpi).
\end{equation*}

Regard $\TypmclEomegazero(=\oplus_{\lambda\in\TypfkAoriginalpi}\bZ\Typeomegazero_\lambda)$ as the unital associative commutative $\bZ$-algebra
defined by $\Typeomegazero_\lambda\Typeomegazero_\mu:=\Typeomegazero_{\lambda+\mu}$,
where the unit is $\Typeomegazero_0$.
Regard $\TyphatmclELam$ as the left $\TypmclEomegazero$-module defined by 
$\Typeomegazero_\lambda\cdot\TypeLam_\mu:=\TypeLam_{\lambda+\mu}$.
By a standard argument, we have
{\begin{equation}\label{eqn:finitekeyeq}
(\prod_{\beta\in\TypRchipipre}(\Typeomegazero_0-\Typeomegazero_{-\beta}))\cdot [\![\nu]\!]^\TypLam
=(\prod_{\al\in\TypRchipipnu}(\Typeomegazero_0+\Typeomegazero_{-\al}
+\cdots +\Typeomegazero_{-(\Typcal-1)\al}))\cdot\TypeLam_\nu
\end{equation}
%For $w\in\TypWchipire$ and $\nu\in\TypfkAoriginalpi$, 
%define $\TypetachipiLam_w(\nu)\in\TypfkAoriginalpi$ 
%(resp. $\TypdotetachipiLam_w(\nu)\in\TypfkAoriginalpi$) by 
%\begin{equation*}
%\quad
%\TypxichipiLam(w)(\TypeLam_\nu)=\TypeLam_{\TypetachipiLam_w(\nu)}\quad
%\mbox{(resp. $\TypdotxichipiLam(w)(\TypeLam_\nu)=\TypeLam_{\TypdotetachipiLam_w(\nu)}$)}.
%\end{equation*}
For $w\in\TypWchipire$, let $[\![w]\!]_\TypLam:=[\![\TypdotetachipiLam(w)(0)]\!]^\TypLam$.

\begin{theorem}\label{theorem:Main}
{\rm{(}}Recall the assumptions \eqref{eqn:assumption-a},
\eqref{eqn:assumption-b} and \eqref{eqn:assTyp}.{\rm{)}}
Let $W:=\TypWchipire$. Then we have
\begin{equation*}
[\TypmclLchipiLam]^\TypLam=\sum_{w\in W}\Typsgnchipi(w)[\![w]\!]_\TypLam,
\end{equation*}
where $[\![w]\!]_\TypLam\ne[\![w^\prime]\!]_\TypLam$
$(w\ne w^\prime)$.
In particular,
\begin{equation*}
\begin{array}{lcl}
\dim\TypmclLchipiLam_\lambda
&=&\sum_{w\in W}\Typsgnchipi(w) 
\dim\TypmclM^{\Typchi,\Typpi}(\TypLam^{+\TypdotetachipiLam(w)(0)})_{\lambda-\TypdotetachipiLam(w)(0)} \\
&=&\sum_{w\in W}\Typsgnchipi(w) 
\dim\TypUmchipi_{\lambda-\TypdotetachipiLam(w)(0)}
\end{array}
\end{equation*}
$(\lambda\in\TypfkAoriginalpi)$,
and $-\TypdotetachipiLam(w)(0)\in\TypfkAoriginalpip\setminus\{0\}$
$(w\ne e)$.
\end{theorem}
{\it{Proof.}} 
For $w\in W$, let $p_w:=|\{w^\prime\in W|[\![w^\prime]\!]_\TypLam=[\![w]\!]_\TypLam\}|$.
Then $p_w=p_e$ for all $w\in W$.
By \eqref{eqn:trueinv}, using a standard argument
such as that for \cite[Proposition and Corollary in Subsection~24.2]{Hum},
we see that there exist $b_w\in{\frac 1 {p_e}}\bZ$ $(w\in W)$
satisfying the following $(b1)$-$(b3)$.
\newline\newline
$(b1)$\quad
$b_w=b_{w^\prime}$ for $w$, $w^\prime\in W$
with $[\![w^\prime]\!]_\TypLam=[\![w]\!]_\TypLam$.
\newline
$(b2)$\quad $[\TypmclLchipiLam]^\TypLam=\sum_{w\in W}b_w[\![w]\!]_\TypLam$.
\newline
$(b3)$\quad For $w\in W$, $b_w=0$ if
$-\TypdotetachipiLam(w)(0)\notin\TypfkAoriginalpip$.
\newline\newline
(As for $(b1)$, keep in mind the set $W/H$ of right class classes,
where $H:=\{w^\prime\in W|[\![w^\prime]\!]_\TypLam=[\![w]\!]_\TypLam\}$
for a fixed $w\in W$.)
By \eqref{eqn:finitekeyeq}, we have
\begin{equation}\label{eqn:prkeyEQ}
\begin{array}{l}
(\prod_{\beta\in\TypRchipipre}(\Typeomegazero_0-\Typeomegazero_{-\beta}))\cdot
[\TypmclLchipiLam]^\TypLam \\
\quad =(\prod_{\al\in\TypRchipipnu}(\Typeomegazero_0+\Typeomegazero_{-\al}
+\cdots +\Typeomegazero_{-(\Typcal-1)\al}))\cdot\sum_{w\in W}b_w
\TypeLam_{\TypdotetachipiLam(w)(0)}.
\end{array}
\end{equation}
By \eqref{eqn:eqmqb}, \eqref{eqn:sbetaLamfm} and \eqref{eqn:prkeyEQ},
we see that $b_{\Typsbeta w}=-b_w$ and 
$\TypdotetachipiLam(\Typsbeta w)(0)
\ne\TypdotetachipiLam(w)(0)$
for $\beta\in\TypRchipipre$ and $w\in W$ with $b_w\ne 0$.
Hence  \eqref{eqn:prkeyEQ} implies that $b_e={\frac 1 {p_e}}$ and
\begin{equation}\label{eqn:Keysgn}
\forall w\in W,\,\, b_w={\frac {\Typsgnchipi(w)} {p_e}}\,\,(\ne 0).
\end{equation}
For $\beta\in\TypRchipipre$, we have $k_\beta\in\bN$ with 
$\TypdotxichipiLam(\Typsbeta)(\TypeLam_0)=\TypeLam_{-k_\beta\beta}$,
i.e., $-k_\beta\beta=\TypdotetachipiLam(\Typsbeta)(0)$.
Hence by Lemma~\ref{lemma:colfcW}~(3)
and \eqref{eqn:sbnemu}, we see that 
$-\TypdotetachipiLam(w)(0)\in\TypfkAoriginalpip\setminus\{0\}$
for all $w\in W\setminus\{e\}$.
Hence $p_e=1$.
\hfill $\Box$

\begin{remark}\label{remark:rankonecase}
Assume that $\Typqbeta=1$ for some 
$\beta\in\TypRchipip$.
Recall that $(\Typchi,\Typpi)$ is irreducible (see \eqref{eqn:assumption-a}).
Then $\Typpibar$ is of $(\Typpibar 0)$, see \eqref{eqn:qbetak}.
We have $\TyphatPchipi=0$.
We see that the following (a), (b), (c) are equivalent.
\begin{equation*}
\mbox{(a) $\dim\TypmclLchipiLam<\infty$.
(b) $\dim\TypmclLchipiLam=1$.
(c) $\TypLam(\TyptrK_{\Typpial_1}\TyptrL_{-\Typpial_1})=1$.}
\end{equation*}
\end{remark}
\subsection{General case} \label{subsection:Gen}
Assume ${\mathrm{Char}}(\bK)=0$. 
Let $\Typchi$ be as in \eqref{eqn:bich}.
Assume $|\TypRchipip|<\infty$.
In this subsection, we do not impose the assumption
\eqref{eqn:assumption-a} or the one \eqref{eqn:assumption-b}.
Nevertheless the notations in this subsection can be defined in the same way as above.
Recall \eqref{eqn:ident}.
Let $\TypUochipi^{\mathrm{ex}}:=\oplus_{\lambda,\mu\in\TypfkAoriginalEX}\bK\TyptrK_\lambda\TyptrL_\mu$.
Let $\Typchi^\prime:=\Typchi_{|\TypfkAoriginalpi\times\TypfkAoriginalpi}$.
Then we have the $\bK$-linear isomorphism
$f_1:\TypU(\Typchi^\prime,\Typpi)\otimes\TypUochipi^{\mathrm{ex}}\to\TypUchipi$
defined by $f_1(X\otimes Y):=XY$.
In particular, $\TypR^{\Typchi^\prime,\pi}_+=\TypRchipip$.
Let $\TypLam:\TypUochipi\to\bK$ be a $\bK$-algebra homomorphism.
Let $\TypLam^\prime:=\TypLam_{|\TypU(\Typchi^\prime,\Typpi)}$.
Then we have the $\bK$-linear isomorphism
$f_2:{\mathcal{L}}^{\Typchi^\prime,\Typpi}(\TypLam^\prime)\to\TypmclLchipiLam$
with $f_2(X v_{\TypLam^\prime})=X v_{\TypLam^\prime}$ $(X\in\TypU(\Typchi^\prime,\Typpi))$. 
In particular, 
\begin{equation} \label{eqn:relEX}
\dim{\mathcal{L}}^{\Typchi^\prime,\Typpi}(\TypLam^\prime)_\lambda
=\dim\TypmclLchipiLam_\lambda\quad(\lambda\in\TypfkAoriginalpi).
\end{equation}
Assume that there exists a non-empty proper subset of $\TypfkI_1$ of $\TypfkI$ such that
$\Typq_{ij}\Typq_{ji}=1$
$(i\in\TypfkI_1, j\in\TypfkI\setminus\TypfkI_1)$.
Let $\TypfkI_2:=\TypfkI\setminus\TypfkI_1$.
For $t\in\TypfkJ_{1,2}$, 
let $\TypfkAoriginal_t:=\oplus_{i\in\TypfkI_t}\bZ\Typpial_i$,
$\Typpi_t:=\Typpi_{\TypfkI_t}$
and $\Typchi_t:=\Typchi_{|\TypfkAoriginal_t\times\TypfkAoriginal_t}$
and $\TypLam_t:=\TypLam_{|\TypU(\Typchi_t,\Typpi_t)}$.
Then we have the $\bK$-linear isomorphisms
\begin{equation*}
\begin{array}{l}
f_3:\TypU(\Typchi_1,\Typpi_1)\otimes\TypU(\Typchi_2,\Typpi_2)\otimes
\TypUochipi^{\mathrm{ex}}\to\TypUchipi \\
\quad (f_3(X_1\otimes X_2\otimes Y):=X_1X_2Y),
\\
f_4:{\mathcal{L}}^{\Typchi_1,\Typpi_1}(\TypLam_1)
\otimes{\mathcal{L}}^{\Typchi_2,\Typpi_2}(\TypLam_2)
\to\TypmclLchipiLam
\quad (f_4(X_1v_{\TypLam_1}\otimes X_2v_{\TypLam_2}):=X_1X_2v_\TypLam),
\end{array}
\end{equation*}
where $X_t\in\TypU(\Typchi_t,\Typpi_t)\,(t\in\TypfkJ_{1,2})$ and $Y\in\TypUochipi^{\mathrm{ex}}$.
In particular, we have
\begin{equation*}
\begin{array}{l}
\TypRchipip=\TypR^{\Typchi_1,\pi_1}_+\cup\TypR^{\Typchi_2,\pi_2}_+,
\,\, \TypR^{\Typchi_1,\pi_1}_+\cap\TypR^{\Typchi_2,\pi_2}_+=\emptyset, \\
\dim\TypmclLchipiLam_{\lambda_1+\lambda_2}
=\dim{\mathcal{L}}^{\Typchi_1,\Typpi_1}(\TypLam_1)_{\lambda_1}
\cdot\dim{\mathcal{L}}^{\Typchi_2,\Typpi_2}(\TypLam_2)_{\lambda_2}$ $(\lambda_t\in\TypfkAoriginal_t\,(t\in\TypfkJ_{1,2}).
\end{array}
\end{equation*}
By \cite[Theorem~7.3]{HY10}, we see that
\begin{equation}\label{eqn:Shapo}
\begin{array}{l}
\mbox{if $\Typqbeta\ne 1$ and $\TypLam(\Typhrhochipibeta\Typqbeta^{-t}\TyptrK_\beta\TyptrL_\beta^{-1}-1)\ne 0$} \\
\mbox{for all $\beta\in\TypRchipip$ and all $t\in\bN$
with $\Typqbeta^t\ne 1$,} \\
\mbox{then $\dim\TypmclLchipiLam_\lambda=\dim\TypmclMchipiLam_\lambda
=\dim\TypUpchipi_\lambda$
for all $\lambda\in\TypfkAoriginalpi$.}
\end{array}
\end{equation}
For $\beta\in\TypRchipipnu$, we
have
$\TypPchipipmbeta=\prod_{t=1}^{\Typcbeta-1}((\Typhrhochipibeta\TyptrK_\beta\TyptrL_{-\beta})^{\pm 1}-
\Typqbeta^t)$.

For $\beta\in\TypRchipire$, we have
the $\bZ$-module automorphism 
${\bar{\Typs}}_\beta:\TypfkAoriginalpi\to\TypfkAoriginalpi$
defined in the way that ${\bar{\Typs}}_\beta(\lambda):=\lambda-k\beta$
for $\lambda\in\TypfkAoriginalpi$ with $k\in\bZ$ satisfying
$\Typchi(\beta,\lambda)\Typchi(\lambda,\beta)=\Typq_\beta^k$.
Let ${\bar{\TypW}}^{\Typchi,\Typpi}_{{\mathrm{real}}}$ be the subgroup of ${\mathrm{Aut}}_\bZ(\TypfkAoriginalpi)$
generated by 
$\{{\bar{\Typs}}_\beta|\beta\in\TypRchipire\}$.
By Lemma~\ref{lemma:colfcW}~(5),
we have the group homomorphism ${\overline{{\mathrm{sgn}}}}^{\Typchi,\Typpi}:
{\bar{\TypW}}^{\Typchi,\Typpi}_{{\mathrm{real}}}\to\{-1,1\}$
such that ${\overline{{\mathrm{sgn}}}}^{\Typchi,\Typpi}({\bar{\Typs}}_\beta)=-1$
$(\beta\in\TypRchipire)$.
For $\beta\in\TypRchipire$,
define $r_\beta(\in\bZ)$ by $\Typhrhochipibeta=\Typqbeta^{r_\beta}$
(see also \eqref{eqn:rbeta}),
and, if $\dim\TypmclLchipiLam<\infty$, 
let $n^\TypLam_\beta
(\in\bZ)$ be such that 
$\TypLam(\TyptrK_\beta\TyptrL_{-\beta})=\Typqbeta^{n^\TypLam_\beta}$.

Using the above argument, by Theorem~\ref{theorem:Main}, we see that
\begin{theorem} \label{theorem:ScMain}
Let $\Typchi$ and $\TypLam$ be of this subsection.
Assume that $\dim\TypmclLchipiLam<\infty$
and $\TypLam(\TyphatPchipi)\ne 0$.
Then we have the group action of ${\bar{\TypW}}^{\Typchi,\Typpi}_{{\mathrm{real}}}$ on
$\TypfkAoriginalpi$ defined by
${\bar{\Typs}}_\beta\cdot\lambda:={\bar{\Typs}}_\beta(\lambda)
-(r_\beta+n^\TypLam_\beta)\beta$
$(\beta\in\TypRchipire,\lambda\in\TypfkAoriginalpi)$. Moreover we
have
\begin{equation*}
\dim\TypmclLchipiLam_\lambda
=\sum_{w\in{\bar{\TypW}}^{\Typchi,\Typpi}_{{\mathrm{real}}}}{\overline{{\mathrm{sgn}}}}^{\Typchi,\Typpi}(w) 
\dim\TypUmchipi_{\lambda-w\cdot 0}
\quad (\lambda\in\TypfkAoriginalpi),
\end{equation*} where
$-w\cdot 0\in\TypfkAoriginalpip\setminus\{0\}$
for $w\ne e$.
\end{theorem}

\section{Appendix---Classification of $\TypLam$ with $\dim\TypmclLchipiLam$ $<\infty$
for $\Typpi=\Typpibar$}
Keep the notation as in Subsection~\ref{subsection:datum}.
In particular, $m$ and $n$ $(\in\bZgeqo)$ mean those of $(\Typpibar 2)$.

Let $\Typpialbar_0\,(\in\TypfkAoriginalbarpip)$ be $0$
(resp. $\sum_{i=n}^\Typl\Typpialbar_i$,
resp. $\Typpialbar_{\Typl-1}+\Typpialbar_\Typl+2\sum_{i=n}^{\Typl-2}\Typpialbar_i$,
resp. $2\Typpialbar_1+3\Typpialbar_2+2\Typpialbar_3+\Typpialbar_4$,
resp. $\Typpialbar_1+2\Typpialbar_2+\Typpialbar_3$,
resp. $\Typpialbar_1+3\Typpialbar_2+2\Typpialbar_3+\Typpialbar_4$,
resp. $\Typpialbar_1+2\Typpialbar_2$)
if $\Typpibar$ is of $(\Typpibar 0)$, $(\Typpibar 1)$, $(\Typpibar 2)$-{\rm{(i)}}, $(\Typpibar 2)$-{\rm{(iii)}} or $(\Typpibar 3)$-{\rm{(ii)}
(resp. $(\Typpibar 2)$-{\rm{(ii)}},
resp. $(\Typpibar 2)$-{\rm{(iv)}}, 
resp. $(\Typpibar 2)$-{\rm{(v)}}
resp. $(\Typpibar 2)$-{\rm{(vi)}}, $(\Typpibar 2)$-{\rm{(vii)}} or $(\Typpibar 3)$-{\rm{(i)},
resp. $(\Typpibar 4)$,
resp. $(\Typpibar 5)$).
Let ${\bar{\Pi}}_0\,(\subset\TypfkAoriginalbarpip)$ be
$\emptyset$
(resp. $\{\Typpialbar_i|i\in\TypfkI\}$,
resp. $\{\Typpialbar_i|i\in\TypfkI\setminus\{m+1\}\}$,
resp. $\{\Typpialbar_0\}\cup\{\Typpialbar_i|i\in\TypfkI\setminus\{n\}\}$,
resp. $\{\Typpialbar_i|i\in\TypfkI\setminus\{1\}\}$,
resp. $\{\Typpialbar_0,\Typpialbar_2,\Typpialbar_3,\Typpialbar_4\}$,
resp. $\{\Typpialbar_0,\Typpialbar_2,\Typpialbar_3\}$,
resp. $\{\Typpialbar_0,\Typpialbar_1,\Typpialbar_3\}$,
resp. $\{\Typpialbar_1,\Typpialbar_3\}$,
resp. $\{\Typpialbar_0,\Typpialbar_1,\Typpialbar_2,\Typpialbar_4\}$,
resp. $\{\Typpialbar_0,\Typpialbar_1\}$)
if $\Typpibar$ is of $(\Typpibar 0)$
(resp. $(\Typpibar 1)$,
resp. $(\Typpibar 2)$-{\rm{(i)}},
resp. $(\Typpibar 2)$-{\rm{(ii)}} or $(\Typpibar 2)$-{\rm{(vi)}},
resp. $(\Typpibar 2)$-{\rm{(iii)}},
resp. $(\Typpibar 2)$-{\rm{(v)}},
resp. $(\Typpibar 2)$-{\rm{(vi)}},
resp. $(\Typpibar 2)$-{\rm{(vii)}} or $(\Typpibar 3)$-{\rm{(i)},
resp. $(\Typpibar 3)$-{\rm{(ii)}, 
resp. $(\Typpibar 4)$,
resp. $(\Typpibar 5)$). 
We see that 
\begin{equation*}
\begin{array}{l}
\mbox{if $\Typpibar$ is not of $(\Typpibar 0)$, then ${\bar \Pi}_0$ is the base of} \\
\mbox{the root system $\TypRchibarpire$
with $\TypRchibarpipre=\TypRchibarpire\cap\TyprmSpan_{\bZgeqo}{\bar \Pi}_0$.}
\end{array}
\end{equation*}
If $\Typpialbar_0=0$, let $c_\Typpibar:=0\,(\in\bZgeqo)$.
If $\Typpialbar_0\ne 0$, then let  $c_\Typpibar\,(\in\bZgeqo)$ be 
$2m$
(resp. $m$, 
resp. $4$, 
resp. $6$, 
resp. $2$, 
resp. $3$)
if $\Typpibar$ is of $(\Typpibar 2)$-{\rm{(ii)}}
(resp. $(\Typpibar 2)$-{\rm{(iv)}}, 
resp. $(\Typpibar 2)$-{\rm{(v)}}, 
resp. $(\Typpibar 2)$-{\rm{(vi)}}
resp. $(\Typpibar 2)$-{\rm{(vii)}}, $(\Typpibar 3)$-{\rm{(i)} or $(\Typpibar 5)$, 
resp. $(\Typpibar 4)$).
\begin{theorem}{\rm{(\cite{AYY15})}}\label{theorem:MainAYY15}
Assume $\Typpibar=\Typpi$.
Let $\TypLam:\TypUochipi\to\bK$ be a $\bK$-algebra homomorphism.
Let $\lambda_\beta:=\TypLam(\TyptrK_\beta\TyptrL_\beta^{-1})$ 
$(\beta\in{\bar{\Pi}}_0\cup\{\Typpialbar_0\})$.
\newline
{\rm{(1)}} If $\dim\TypmclLchipiLam<\infty$, then
\begin{equation}\label{eqn:cdint}
\forall\beta\in{\bar{\Pi}}_0\cup\{\Typpialbar_0\}, \exists t_\beta\in\bZgeqo, \lambda_\beta=\Typqbeta^{t_\beta}.
\end{equation}
\newline
{\rm{(2)}} Assume that $\Lambda$ satisfies the condition in \eqref{eqn:cdint}.
Assume $t_{\Typpialbar_0}=0$ if $\Typpialbar_0=0$. 
Then 
\begin{equation*}
\mbox{$\dim\TypmclLchipiLam<\infty$ if $t_{\Typpialbar_0}\geq c_\Typpibar$.}
\end{equation*}
Moreover, if $t_{\Typpialbar_0}<c_\Typpibar$,
then $\dim\TypmclLchipiLam<\infty$ if and only if one of the following
$(C1)$-$(C11)$ holds. \newline\newline
$(C1)$ $\Typpibar$ is of $(\Typpibar 2)$-{\rm{(ii)}}. 
$t_{\Typpialbar_0}\in\TypfkJ_{0,2m-2}\cap 2\bZ$.
$\lambda_{{\Typpialbar_i}}=1$ for all $i\in\TypfkJ_{n+{\frac {t_{\Typpialbar_0}} 2}+1,\Typl}$.
\newline
$(C2)$ $\Typpibar$ is of $(\Typpibar 2)$-{\rm{(iv)}}. 
$t_{\Typpialbar_0}\in\TypfkJ_{0,m-2}$.
$\prod_{i=n}^{n+t_{\Typpialbar_0}}\lambda_{\Typpialbar_i}=q^{-2t_{\Typpialbar_0}}$.
$\lambda_{{\Typpialbar_j}}=1$ for all $j\in\TypfkJ_{n+t_{\Typpialbar_0}+1,\Typl}$.
\newline
$(C3)$ $\Typpibar$ is of $(\Typpibar 2)$-{\rm{(iv)}}. 
$t_{\Typpialbar_0}=m-1$.
$\prod_{i=n}^{\Typl-1}\lambda_{\Typpialbar_i}=q^{-2(m-1)}$.
\newline
$(C4)$ $t_{\Typpialbar_0}=0$. $\lambda_{{\Typpialbar_i}}=1$ for all $i\in\TypfkI$.
\newline
$(C5)$ $\Typpibar$ is of $(\Typpibar 2)$-{\rm{(v)}}. $t_{\Typpialbar_0}=2$.
$\lambda_{{\Typpialbar_2}}=\lambda_{{\Typpialbar_4}}=1$.
$\lambda_{{\Typpialbar_1}}\lambda_{{\Typpialbar_3}}=q^{-6}$.
\newline
$(C6)$ $\Typpibar$ is of $(\Typpibar 2)$-{\rm{(v)}}. $t_{\Typpialbar_0}=3$.
$\lambda_{{\Typpialbar_1}}\lambda_{{\Typpialbar_3}}\lambda_{{\Typpialbar_4}}^2=q^{-12}$.
$\lambda_{{\Typpialbar_2}}=q^2\lambda_{{\Typpialbar_3}}$.
\newline
$(C7)$ $\Typpibar$ is of $(\Typpibar 2)$-{\rm{(vi)}}. 
$t_{\Typpialbar_0}=4$. $\lambda_{{\Typpialbar_2}}=1$.
\newline
$(C8)$ $\Typpibar$ is of $(\Typpibar 2)$-{\rm{(vii)}}. 
$t_{\Typpialbar_0}=1$. $\lambda_{{\Typpialbar_2}}=1$ or $\lambda_{{\Typpialbar_1}}\lambda_{{\Typpialbar_2}}=q^2$. 
\newline
$(C9)$ $\Typpibar$ is of $(\Typpibar 3)$-{\rm{(i)}}. 
$t_{\Typpialbar_0}=1$. $\lambda_{{\Typpialbar_2}}=1$ or $\lambda_{{\Typpialbar_1}}\lambda_{{\Typpialbar_2}}=(xy)^{-1}$. 
\newline
$(C10)$ $\Typpibar$ is of $(\Typpibar 4)$. 
$t_{\Typpialbar_0}=1$. $\lambda_{{\Typpialbar_3}}=\lambda_{{\Typpialbar_3}}=1$. 
\newline
$(C11)$ $\Typpibar$ is of $(\Typpibar 4)$. 
$t_{\Typpialbar_0}=2$. $\lambda_{{\Typpialbar_3}}=1$ or 
$\lambda_{{\Typpialbar_1}}\lambda_{{\Typpialbar_3}}\lambda_{{\Typpialbar_4}}=1$. 
\end{theorem}

\begin{remark}\label{remark:Geer}
We make a comment on the N.~Geer's result \cite{Geer07}.
Keep the notation of Introduction.
Recall $U_q({\mathfrak{g}})$ and $\Omega\in{\mathrm{Gr}}_\Pi$ from Introduction.
Assume that $\dim{\mathfrak{g}}<\infty$ and
$(\al_i,\al_j)\in\bZ$ $(i,j\in\TypfkI)$.
Then $q\in\bCtinf$.
If 
$q$ is transcendental over $\bQ$ and the condition is satisfied that
\begin{equation}\label{eqn:GeerCond}
\forall i\in\TypfkI, \exists a_i\in\bZ,\,\,\mbox{s.t.}\,\,\Omega(\al_i)=\Typq^{a_i},
\end{equation}
then \cite[Theorem~1.2]{Geer07}
tells that 
the dimensions $\dim\TypIntmcLOm_\lambda$
$(\lambda\in\bZgeqo\Pi)$
are the same as those for
the irreducible highest weight module of ${\mathfrak{g}}$
having the highest weight vector $v$
with $h_iv=a_iv$ $(i\in\TypfkI)$.
However, it follows from
Theorem~\ref{theorem:MainAYY15} that there are
$\Omega$ with $\dim\TypIntmcLOm<\infty$
such that $\Omega$ do not satisfy the condition \eqref{eqn:GeerCond}
(i.e., there cannot exist 
an irreducible highest weight module of ${\mathfrak{g}}$ obtained from such $\Omega$ 
by taking $q\rightarrow 1$).
For example, we have the following {\rm{(1)}}-{\rm{(3)}}. 
Let $\Typpibar$ be the one of Theorem~\ref{theorem:MainAYY15}
corresponding to ${\mathfrak{g}}$. 
Assume that $\Typpibar$ is of $(\Typpibar 2)$-{\rm{(t)}}
and $(\al_i,\al_j)$ are $x_{ij}$ of $(\Typpibar 2)$-{\rm{(t)}},
where $(t)$ is one of ${\rm{(i)}},\ldots,{\rm{(vii)}}$.
Let $\Omega_i:=\Omega(2\al_i)$ $(i\in\TypfkI)$.
Let $g_\Omega:=\prod_{\gamma\in\TypIntRpnull}(1-(-1)^{p(\gamma)}{\hat{\rho}}_\Pi(\gamma)\Omega(2\gamma))$,
where recall \eqref{eqn:AssScMain}.
Identify $\al_i$ $(i\in\TypfkI)$ of Introduction with $\Typpialbar_i=\Typpibar(i)$ of Subsection~\ref{subsection:datum}.
Recall the Shapovalov determinant formula \cite[Theorem~7.3]{HY10} and see also \cite[(7.14) and Theorem~7.5]{BY18}.
\newline
{\rm{(1)}} Assume that $(t)$ is {\rm{(ii)}}.
Then ${\mathfrak{g}}$ is isomorphic to $B(m,n)=osp(2m+1|2n)$
$(m,n\in\bN)$,
and we have $\TypIntRpreal=\{\Typepsilon_{i_1}\,(i_1\in\TypfkJ_{1,m+n}),\,
\Typepsilon_{i_2}\pm\Typepsilon_{j_2}\,(\{i_2,j_2\}\subset\TypfkJ_{1,n},i_2<j_2,
\,\mbox{or}\,\{i_2,j_2\}\subset\TypfkJ_{n+1,m+n},i_2<j_2)\}$
and 
$\TypIntRpnull=\{\Typepsilon_{i_3}\pm\Typepsilon_{j_3}|
i_3\in\TypfkJ_{1,n},\,i_3\in\TypfkJ_{n+1,m+n}\}$,
where $\Typepsilon_i:=\sum_{k=i}^{m+n}\al_k$ $(i\in\TypfkI)$.
In particular, $\{{\hat{\rho}}_\Pi(\beta)|\beta\in\TypIntRpnull\}
\subset\{-\Typq^{2b}|b\in\bZ\}$.
If there 
exist $d_i\in\bZgeqo$ $(i\in\TypfkI)$
with 
$\Omega_y=\Typq^{-2d_y}$ $(y\in\TypfkJ_{1,n-1})$,
$\Omega_z=\Typq^{2d_z}$ $(z\in\TypfkJ_{n+1,m+n-1})$,
$\Omega_{m+n}=\Typq^{d_{m+n}}$
and $\prod_{j=n}^{m+n}\Omega_j=-\Typq^{-(2(m+d_n)+1)}$, then 
$\Omega$ does not satisfy the condition \eqref{eqn:GeerCond}
but Theorem~\ref{theorem:MainAYY15}~(2) tells $\dim\TypIntmcLOm<\infty$ 
and
$g_\Omega\ne 0$,
which implies that the equation \eqref{eqn:FormulaScMain} holds.
\newline
{\rm{(2)}} Assume that $(t)$ is {\rm{(vii)}}.
Then ${\mathfrak{g}}$ is isomorphic to $D(2,1;a)$ with $a\in\bZ\setminus\{0,-1\}$,
where $D(2,1;1)$ is isomorphic to $D(2,1)$.
Let $\beta_1:=\al_1$, $\beta_2=\al_1+2\al_2+\al_3$, $\beta_3:=\al_3$
and $\beta_4:=\al_2$,  $\beta_5:=\al_1+\al_2$, $\beta_6:=\al_2+\al_3$, $\beta_7:=\al_1+\al_2+\al_3$.
Then $\TypIntRpreal=\{\beta_1,\beta_2,\beta_3\}$ and 
$\TypIntRpnull=\{\beta_4,\beta_5,\beta_6,\beta_7\}$.
Note $-2(\al_2,\al_i)=(\al_i,\al_i)$ $(i\in\TypfkI)$
(recall from \cite[Subsection~8.3 and Corollary~8.5.4]{Musson12} that the $\rho$ of a basic Lie superalgebra is ${\frac 1 2}(\sum_{\beta\in\TypIntRpreal}\beta-
\sum_{\gamma\in\TypIntRpnull}\gamma)$).
Then
$g_\Omega
=\prod_{\gamma\in\TypIntRpnull}(1+{\hat{\rho}}_\Pi(\gamma)\Omega(2\gamma))
=\prod_{j=4}^7(1-q^{-2(\al_2,\beta_j)}\Omega(2\beta_j))
=(1-\Omega_2)(1-q^{-2}\Omega_1\Omega_2)(1-q^{-2a}\Omega_2\Omega_3)
(1-q^{-2(1+a)}\Omega_1\Omega_2\Omega_3)$.
Assume that $\Omega_1=q^{-2d_1}$, $\Omega_3=q^{-2ad_3}$
and $\Omega_1\Omega_2^2\Omega_3=q^{2(1+a)(2+d_2)}$ 
for some $d_t\in\bZgeqo$ $(t\in\TypfkJ_{1,3}$).
Then $\dim\TypIntmcLOm<\infty$ by Theorem~\ref{theorem:MainAYY15}~(2).
We have $\Omega_2=(-1)^uq^{d_1+ad_3+(1+a)(2+d_2)}$
for some $u\in\TypfkJ_{0,1}$.
Then $g_\Omega=0$ if and only if
$u=0$ and 
$0\in\{d_1+ad_3+(1+a)(2+d_2),
-2-d_1+ad_3+(1+a)(2+d_2),
-2a+d_1-ad_3+(1+a)(2+d_2),
-d_1-ad_3+(1+a)d_2\}$.
Then \cite[Theorem~7.3]{HY10} tells that
for $\lambda\in\bZgeqo\Pi$, $\dim\TypIntmcLOm_\lambda<\dim\TypIntmcMOm_\lambda$
if and only if there exists $((b_1,b_2,b_3),(b_4,b_5,b_6,b_7))\in(\bZgeqo)^3\times(\TypfkJ_{0,1})^4$
such that $\sum_{k=1}^7b_k\beta_k=\lambda$ and
$\prod_{k=1}^7\prod_{t_k=1}^{b_k}(q^{t_k(\beta_k,\beta_k)}-q^{-2(\al_2,\beta_k)}\Omega(2\beta_k))
= 0$. Note that
the LHS of the last equation equals
\begin{equation*}
\begin{array}{l}
(\prod_{t_1=1}^{b_1}(q^{-2t_1}-q^{-2(1+d_1)}))
(\prod_{t_2=1}^{b_2}(q^{2(1+a)t_2}-q^{2(1+a)(1+d_2)})) \\
\cdot(\prod_{t_3=1}^{b_3}(q^{-2at_3}-q^{-2a(1+d_3)})) 
(\prod_{t_4=1}^{b_4}(1-(-1)^uq^{d_1+ad_3+(1+a)(2+d_2)})) \\
\cdot(\prod_{t_5=1}^{b_5}(1-(-1)^uq^{-(2+d_1)+ad_3+(1+a)(2+d_2)})) \\
\cdot(\prod_{t_6=1}^{b_6}(1-(-1)^uq^{d_1-a(2+d_3)+(1+a)(2+d_2)})) 
(\prod_{t_7=1}^{b_7}(1-(-1)^uq^{-d_1-ad_3+(1+a)d_2})).
\end{array}
\end{equation*}
Let $k\in\bN$.
Assume that $d_1=d_2=d_3=k$ and $\lambda=k(\al_1+\al_2+\al_3)$.
If $u=0$ (resp. $u=1$),
then $g_\Omega=0$ and $\dim\TypIntmcLOm_\lambda<\dim\TypIntmcMOm_\lambda$
(resp. $g_\Omega\ne 0$ and $\dim\TypIntmcLOm_\lambda =\dim\TypIntmcMOm_\lambda$).
\newline
{\rm{(3)}} Assume that $(t)$ is {\rm{(iv)}}.
Then ${\mathfrak{g}}$ is isomorphic to $D(m,n)=osp(2m|2n)$
$(m\in\TypfkJ_{2,\infty},n\in\bN)$,
and we have $\TypIntRpreal=\{2\Typepsilon_{i_1}\,(i_1\in\TypfkJ_{1,n}),\,
\Typepsilon_{i_2}\pm\Typepsilon_{j_2}\,(\{i_2,j_2\}\subset\TypfkJ_{1,n},i_2<j_2,
\,\mbox{or}\,\{i_2,j_2\}\subset\TypfkJ_{n+1,m+n},i_2<j_2)\}$
and 
$\TypIntRpnull=\{\Typepsilon_{i_3}\pm\Typepsilon_{j_3}|
i_3\in\TypfkJ_{1,n},\,i_3\in\TypfkJ_{n+1,m+n}\}$,
where $\Typepsilon_i:={\frac 1 2}(\al_{m+n}-\al_{m+n-1})+\sum_{k=i}^{m+n-1}\al_k$ $(i\in\TypfkI)$.
In particular, $\{{\hat{\rho}}_\Pi(\beta)|\beta\in\TypIntRpnull\}
\subset\{-\Typq^{2b}|b\in\bZ\}$.
Let $\Omega$ be such that there exist $d_i\in\bZgeqo$ $(i\in\TypfkI)$
and $u\in\TypfkJ_{0,1}$ 
with 
$\Omega_y=\Typq^{-2d_y}$ $(y\in\TypfkJ_{1,n-1})$,
$\Omega_z=\Typq^{2d_z}$ $(z\in\TypfkJ_{n+1,m+n})$,
and $\prod_{j=n}^{m+n-1}\Omega_j=\Omega_{m+n}\cdot\prod_{j=n}^{m+n-2}\Omega_j=(-1)^u
\Typq^{-2(m+d_n)}$.
Then $\dim\TypIntmcLOm<\infty$ by Theorem~\ref{theorem:MainAYY15}~(2).
If $u=0$ (resp. $u=1$), $\Omega$ satisfies (resp. does not satisfy) the condition \eqref{eqn:GeerCond}.
If $u=1$, then  $g_\Omega\ne 0$,
which implies that the equation \eqref{eqn:FormulaScMain} holds.
Let $k\in\bN$.
Assume that $d_j=k$ for all $j\in\TypfkJ_{n,m+n}$.
Then $\Omega_n=(-1)^uq^{-2m(k+1)}$.
Let $\lambda:=k\cdot\sum_{j=n}^{m+n}\al_j\,(\in\bZ\Pi)$.
By \cite[Theorem~7.3]{HY10}, we see that 
if $u=0$ (resp. $u=1$),
then $g_\Omega=0$ and $\dim\TypIntmcLOm_\lambda<\dim\TypIntmcMOm_\lambda$
(resp. $g_\Omega\ne 0$ and $\dim\TypIntmcLOm_\lambda =\dim\TypIntmcMOm_\lambda$).
\end{remark}

\vspace{2cm}

Hiroyuki Yamane \par 
Department of Mathematics \par 
Faculty of Science \par  
Toyama University \par  
Gofuku, Toyama 930-8555, JAPAN \par 
e-mail: hiroyuki@sci.u-toyama.ac.jp

\end{document}